\documentclass{amsart}
\usepackage{amsmath,amssymb,amsxtra,amsthm,amscd}
\usepackage[mathscr]{eucal}
\usepackage{bbm}
\usepackage{graphicx}
\usepackage[all,cmtip]{xy}
\usepackage{lmodern}
\usepackage[T1]{fontenc}
\usepackage{microtype}
\usepackage{adjustbox}
\usepackage{color}

\newif\ifShowLabels
\ShowLabelstrue
\newcommand{\TeXref}[1]{
\marginpar{\scriptsize \texttt{#1}}}

\ShowLabelsfalse

\DeclareMathOperator{\A}{\mathbf{A}}

\DeclareMathOperator{\B}{\mathbf{B}}

         \newcommand{\BX}{\B_{X}}

\DeclareMathOperator{\bsv}{\boldsymbol{v}}

\DeclareMathOperator{\C}{\mathbf{C}}

\DeclareMathOperator{\co}{\mathrm{co}}
\DeclareMathOperator{\coim}{coim}
\DeclareMathOperator{\coker}{coker}

\DeclareMathOperator{\D}{\mathbf{D}}

\DeclareMathOperator{\Dw}{\mathbf{D}^{\boldsymbol{w}}}
\DeclareMathOperator{\E}{\mathbf{E}}

\DeclareMathOperator{\fil}{fil}

\DeclareMathOperator{\Free}{\mathbf{Free}}

\DeclareMathOperator{\G}{\mathbf{G}}

         \newcommand{\Gnc}{G^{-\infty}}

         \newcommand{\GncX}{{G}_{X}^{-\infty}}

\DeclareMathOperator{\Hom}{Hom}

\DeclareMathOperator{\id}{id}

\DeclareMathOperator{\im}{im}

\DeclareMathOperator{\K}{\mathit{K}}
         \newcommand{\Knc}{\K^{-\infty}}

\DeclareMathOperator{\LI}{\mathbf{L}}
\DeclareMathOperator{\SI}{\mathbf{S}}

\DeclareMathOperator{\BLI}{\mathbf{BL}}
\DeclareMathOperator{\BSI}{\mathbf{BS}}

\DeclareMathOperator{\Mod}{\mathbf{Mod}}

\DeclareMathOperator{\point}{pt}

\DeclareMathOperator{\pr}{pr}

         \newcommand{\subdot}{\boldsymbol{\cdot}}

\DeclareMathOperator{\U}{\mathbf{U}}

\DeclareMathOperator{\vG}{\mathbf{\boldsymbol{v}G}}
\DeclareMathOperator{\wG}{\mathbf{\boldsymbol{w}G}}

\DeclareMathOperator{\w}{\boldsymbol{w}}
\DeclareMathOperator{\vD}{\mathbf{\boldsymbol{v}D}}
\DeclareMathOperator{\vDw}{\mathbf{\boldsymbol{v}D^{\boldsymbol{w}}}}
\DeclareMathOperator{\vE}{\mathbf{\boldsymbol{v}E}}

         \newcommand{\bfw}{\boldsymbol{w}}

\DeclareMathOperator{\wD}{\mathbf{\boldsymbol{w}D}}

\DeclareMathOperator{\Y}{\mathbf{Y}}
\DeclareMathOperator{\Z}{\mathbf{Z}}

\DeclareMathOperator*{\one}{1}
\newcommand{\onehatplace}[1]
{ \one^{\substack{#1 \\ \frown}} }

\DeclareMathOperator*{\bones}{\times}
\newcommand{\undertimes}[1]
{ \bones_{#1} }

\DeclareMathOperator*{\bowl}{\cup}
\newcommand{\undercup}[1]
{ \bowl_{#1} }

\DeclareMathOperator*{\arch}{\cap}
\newcommand{\undercap}[1]
{ \arch_{#1} }

\newcommand{\pull}
{\!\!\! -\!\!\! -\!\!\! -\!\!\!}

\DeclareMathOperator*{\holimprep}{holim}                       
\newcommand{\holim}[1]%
{\displaystyle\holimprep_{\substack{\leftarrow \pull - \\ #1}} \, }

\DeclareMathOperator*{\hocolimprep}{hocolim}                   
\newcommand{\hocolim}[1]%
{\displaystyle\hocolimprep_{\substack{- \pull \rightarrow \\ #1}} \, }

\DeclareMathOperator*{\plainlim}{lim}                           
\newcommand{\contralim}[1]%
{\displaystyle\plainlim_{\substack{\leftarrow \pull - \\ #1}} \, }

\DeclareMathOperator*{\plaincolim}{colim}                       
\newcommand{\colim}[1]%
{\displaystyle\plaincolim_{\substack{- \pull \rightarrow \\ #1}} \, }

\DeclareMathOperator*{\laxlimplain}{laxlim}                     
\newcommand{\laxlim}[1]%
{\displaystyle\laxlimplain_{\substack{\leftarrow \pull - \\ #1}} \, }

\providecommand{\bysame}{\makebox[3em]{\hrulefill}\thinspace}

\swapnumbers
\theoremstyle{plain}
\newtheorem{Thm}{Theorem}[subsection]

\newtheorem{Cor}[Thm]{Corollary}

\newtheorem{Lem}[Thm]{Lemma}
\newtheorem{Prop}[Thm]{Proposition}

\theoremstyle{definition}
\newtheorem{Def}[Thm]{Definition}

\newtheorem{Ex}[Thm]{Example}

\newtheorem{Rem}[Thm]{Remark}

\theoremstyle{remark}
\newtheorem{Not}[Thm]{Notation}

\newtheoremstyle{freestylethm}{6pt}{6pt}{\itshape}{}%
                {\bfseries}{}{.5em}{\thmnote{#3}}
\theoremstyle{freestylethm}

\newcommand{\SecRef}[2]{\section{#1}\label{S:#2}%
\ifShowLabels \TeXref{{S:#2}} \fi}
\newcommand{\SSecRef}[2]{\subsection{#1}\label{SS:#2}%
\ifShowLabels \TeXref{{SS:#2}} \fi}

\newcommand{\refSS}[1]{\textup{\ref{SS:#1}}}

\newcommand{\refT}[1]{\textup{\ref{T:#1}}}
\newcommand{\refL}[1]{\textup{\ref{L:#1}}}
\newcommand{\refD}[1]{\textup{\ref{D:#1}}}
\newcommand{\refC}[1]{\textup{\ref{C:#1}}}
\newcommand{\refE}[1]{\textup{\ref{E:#1}}}
\newcommand{\refP}[1]{\textup{\ref{P:#1}}}
\newcommand{\refR}[1]{\textup{\ref{R:#1}}}

\newenvironment{ThmRef}[1]%
{ \begin{Thm} \label{T:#1}
\ifShowLabels \TeXref{T:#1} \fi }%
{ \end{Thm} }
\newenvironment{DefRef}[1]%
{ \begin{Def} \label{D:#1}
\ifShowLabels \TeXref{D:#1} \fi }%
{ \end{Def} }
\newenvironment{LemRef}[1]%
{ \begin{Lem} \label{L:#1}
\ifShowLabels \TeXref{L:#1} \fi }%
{ \end{Lem} }
\newenvironment{CorRef}[1]%
{ \begin{Cor} \label{C:#1}
\ifShowLabels \TeXref{C:#1} \fi }%
{ \end{Cor} }
\newenvironment{RemRef}[1]%
{ \begin{Rem} \label{R:#1}
\ifShowLabels \TeXref{R:#1} \fi }%
{ \end{Rem} }
\newenvironment{PropRef}[1]%
{ \begin{Prop} \label{P:#1}
\ifShowLabels \TeXref{P:#1} \fi }%
{ \end{Prop} }
\newenvironment{ExRef}[1]%
{ \begin{Ex} \label{E:#1}
\ifShowLabels \TeXref{E:#1} \fi  }%
{ \end{Ex} }
\newenvironment{NotRef}[1]%
{ \begin{Not} \label{N:#1}
\ifShowLabels \TeXref{N:#1} \fi }%
{ \end{Not} }

\newenvironment{ThmRefName}[2]%
{ \begin{Thm} [#2]\label{T:#1}
\ifShowLabels \TeXref{T:#1} \fi }%
{ \end{Thm} }
{ \begin{Def} [#2]\label{D:#1}
\ifShowLabels \TeXref{D:#1} \fi }%
{ \end{Def} }
{ \begin{Lem} [#2]\label{L:#1}
\ifShowLabels \TeXref{L:#1} \fi }%
{ \end{Lem} }
{ \begin{Cor} [#2]\label{C:#1}
\ifShowLabels \TeXref{C:#1} \fi }%
{ \end{Cor} }
{ \begin{Rem} [#2]\label{R:#1}
\ifShowLabels \TeXref{R:#1} \fi }%
{ \end{Rem} }
{ \begin{Prop} [#2]\label{P:#1}
\ifShowLabels \TeXref{P:#1} \fi }%
{ \end{Prop} }
{ \begin{Ex} [#2]\label{E:#1}
\ifShowLabels \TeXref{E:#1} \fi }%
{ \end{Ex} }

\let\oldtocsection=\tocsection
\let\oldtocsubsection=\tocsubsection
\renewcommand{\tocsection}[2]{\hspace{0em}\oldtocsection{#1}{#2}}
\renewcommand{\tocsubsection}[2]{\hspace{2em}\oldtocsubsection{#1}{#2}}

\begin{document}

\title[\textit{G}-theory with fibred control]{Bounded \textit{G}-theory with fibred control}
\author[Gunnar Carlsson]{Gunnar Carlsson}
\address{Department of Mathematics\\ Stanford University\\ Stanford\\ CA 94305}
\email{gunnar@math.stanford.edu}
\author[Boris Goldfarb]{Boris Goldfarb}
\address{Department of Mathematics and Statistics\\ SUNY\\ Albany\\ NY 12222}
\email{goldfarb@math.albany.edu}
\date{\today}

\begin{abstract}
We use filtered modules over a Noetherian ring and fibred bounded control on homomorphisms to construct a new kind of controlled algebra with applications in geometric topology.  The resulting theory can be thought of as a ``pushout'' of bounded $K$-theory with fibred control and bounded $G$-theory constructed and used by the authors.  Bounded $G$-theory was geared toward constructing a $G$-theoretic version of assembly maps and proving the Novikov injectivity conjecture for them.  The $G$-theory with fibred control is needed in the study of surjectivity of the assembly map.  The relation between the $K$- and $G$-theories is the classical one: $K$-theory is meaningful, however $G$-theory is easier to compute, and the relationship is expressed via a Cartan map. This map turns out to be an equivalence under very mild constraints in terms of metric geometry such as finite decomposition complexity.  The fibred theory is certainly more complicated than the absolute theory.  This paper contains the non-equivariant theory including fibred controlled excision theorems known to be crucial for computations. 
\end{abstract}

\maketitle

\tableofcontents

\SecRef{Introduction}{intro}

The purpose of this paper is to use filtered modules over a Noetherian ring with a fibred bounded control on homomorphisms to construct a bounded $G$-theory with fibred control. This theory can be thought of as a ``pushout'' of the bounded $K$-theory with fibred control constructed by the authors in \cite{gCbG:18} and the controlled $G$-theory constructed in \cite{gCbG:00}.  
Here is a summary of this situation:

\begin{figure}[h!]
      \centering
       \[
\xymatrix{
 K (X,R) \ar[r] \ar[d]
&K_X (Y) \ar@{..>}[d] \\
 G (X,R) \ar@{..>}[r]
&{G_X(Y)} }
\]
       \caption{Bounded \textit{G}-theory with fibred control as a ``pushout''.}
      \label{Reed}
   \end{figure}

Throughout this paper, metric spaces such as $X$ and $Y$ that appear in the square will be proper metric spaces in the sense that every closed bounded subspace is compact.  The space or spectrum in the upper left corner represents the indispensable in modern geometric topology bounded $K$-theory of Pedersen and Weibel \cite{ePcW:85,ePcW:89}, reviewed here in the beginning of section \refSS{EBKT}.  This theory is defined for any ring of coefficients $R$.  It is built out of free $R$-modules with generating sets parametrized over the metric space $X$.  This allows to impose geometric control conditions on the homomorphisms $f \colon F \to G$.  The bounded control condition postulates that there is number $b \ge 0$ so that the image of every basis element in $F$ associated to some point $x$ in $X$ is spanned by basis elements in $G$ that are referenced by points within $b$ from $x$.  We will review precise definitions shortly.

The spectrum in the upper right corner $K_X (Y)$ is a generalization of this theory to the situation when the modules are parametrized by the product of two metric spaces $X$ and $Y$, and the control imposed on the homomorphisms is relaxed: it is essentially the bounded control across $X$ but the bound is allowed to change in the complementary direction $Y$ as one varies the $X$-coordinate.  This theory becomes useful when one considers ``bundle phenomena''.   For example, the space $X$ can be the universal cover of the tangent bundle of a manifold embedded in a Euclidean space or even its discrete model such as the fundamental group with a word metric.  The space $Y$ can be the universal cover of the normal bundle to the embedding with a variety of useful metrics.  This situation comes up the authors' work in geometric topology.  The fibred $K$-theory is still defined for any coefficient ring $R$.

To describe the bottom row in the square and for the rest of the paper, we restrict to Noetherian rings $R$.

In place of parametrizations used to control homomorphisms between free modules, one can use filtrations of arbitrary $R$-modules by subsets of the metric space $X$ and impose control conditions in terms of the filtrations.  This was done in \cite{gCbG:00} for a single space $X$.  The result was the \textit{bounded $G$-theory} spectrum $G (X,R)$.  The definition involved promoting the setting from the additive structure for free modules in the definition of bounded $K$-theory to a specific non-split Quillen exact structure on a category of filtered $R$ modules with morphisms satisfying control conditions and the admissible morphisms satisfying further ``bi-control'' conditions.  Regardless of the significant change in techniques, literally every theorem about bounded $K$-theory has an exact (accidental pun) counterpart in $G$-theory.

Now it is clear what the ``pushout'' $G_X(Y)$ is supposed to mean.  We want to look at the $K$-theory of a category built out of filtered modules over the product $X \times Y$ where the morphisms have the fibred control condition of the type described for fibred $K$-theory.  This time we are interested in very specific excision results designed to deconstruct only the ``fiber'' direction.  In our applications of this material we want to perform what we call here \textit{relative excision} in the normal bundle direction.  This is greatly facilitated by the hybrid conditions imposed on the objects themselves.  We include several remarks in the paper regarding the options and why our choices seem to be optimal.  Long story short, we have resolved in this paper the problems that may be much harder to solve, if solvable at all, for the straightforward combination of the theories in the corners of the diagram.  We resolve them for a carefully crafted theory that has all the desired properties and yet specializes to precisely $G (X,R)$ when localized near the subspace $X \times 0$ in $X \times Y$.  

\medskip

\textbf{Acknowledgement.}
We would like to thank the referee for excellent comments that improved the narrative and the precision of the paper.

\SecRef{Elements of bounded \textit{G}-theory}{MSOGS1}

Bounded $G$-theory defined in \cite{gCbG:00} is a variant of bounded $K$-theory of Pedersen and Weibel made applicable to more general, non-split exact  structures.
It was designed by the authors for a different purpose than the one in this paper.  The old focus was on the equivariant theory in addition to very basic excision that were sufficient for introducing an assembly map in $G$-theory for all finitely generated groups and proving the injectivity Novikov type theorem for a large class of groups.
We will review and augment some material from \cite{gCbG:00} in the form best fit for the fibred theory.

\SSecRef{Basic definitions}{EBKT}

We start with a brief recollection of the bounded $K$-theory setup.
The coefficients $R$ for this theory can be an arbitrary associative ring.
The \textit{bounded category} $\mathcal{C} (X,R)$ is the additive category of \textit{geometric $R$-modules} whose 
objects are functions $F \colon X \to \Free_{fg} (R)$ which are locally finite assignments of free finitely generated $R$-modules $F_x$ to points $x$ of $X$.
The local finiteness condition requires precisely that for any bounded subset $S \subset X$ the restriction of $F$ to $X$ has finitely many nonzero values.
Let $d$ be the distance function in $X$.  The morphisms in $\mathcal{C} (M,R)$ are the $R$-linear homomorphisms
\[
\phi \colon \bigoplus_{x \in X} F_x \longrightarrow \bigoplus_{y \in X} G_y
\]
with the property
that the components $F_x \to G_y$ are zero for $d(x,y) > b$
for some fixed real number $b = b (\phi) \ge 0$.
The associated $K$-theory spectrum is denoted by $K (X,R)$
and is called the \textit{bounded K-theory} of $X$.

\begin{RemRef}{HAVRCF}
We would like to remind the reader that the original paper of Pedersen and Weibel \cite{ePcW:85} was already written in greater generality.  If $\mathcal{A}$ is any additive category, it can be used as coefficients in this construction in place of finitely generated free $R$-modules using the same formulas as above.  The outcome is the bounded category $\mathcal{C} (X,\mathcal{A})$ which is again an additive category with the evident notion of split exact sequences.  It is now possible to iterate this construction: when there are two metric spaces $X$ and $Y$, Pedersen and Weibel built the additive category $\mathcal{C} (X,\mathcal{C} (Y,R))$.  The objects of this category can be identified with the objects of $\mathcal{C} (X \times Y,R)$ where the product is given a reasonable metric such as the max metric.  The morphisms are, however, very different.  They are $R$-homomorphisms which are still controlled over $X$ in the standard fashion but the non-zero components $F_x \to G_y$ are now allowed to vary in range with the $X$-coordinates of $x$ and $y$.
In contrast, morphisms in $\mathcal{C} (X \times Y,R)$ require one single number to work as a bound for all of these components.  

Pedersen and Weibel were more interested in the product situation and, in fact, repaired the non-uniform boundedness properties of morphisms in $\mathcal{C} (X,\mathcal{C} (Y,R))$ by filtering the morphism sets.  With that fix the category becomes isomorphic to $\mathcal{C} (X \times Y,R)$.  We, instead, embraced the flexibility of this construction in \cite{gCbG:18} with the idea of exploiting the additional deformations in $K (X,\mathcal{C} (Y,R))$, the \textit{K-theory with fibred control}, that the construction allows.  This is the spectrum that shows up in Figure \ref{Reed} as $K_X (Y)$.
\end{RemRef}

\begin{NotRef}{ttyh}
For a subset $S \subset X$ and a real number $r \ge 0$,
$S[r]$ will stand for the metric $r$-enlargement $\{ x \in X \mid d (x,S) \le r \}$.
In this notation, the metric ball of radius $r$ centered at $x$ is $\{ x \} [r]$ or simply $x [r]$.
\end{NotRef}

A variation of the basic construction of bounded $K$-theory is based on the following observation. 
For every object $F$ and a subset $S$ there is a free $R$-module $F(S) = \bigoplus_{m \in S} F_m$.  In this context we say an
element $x \in F$ is \textit{supported} on a subset $S$ if $x \in F(S)$.
Now the restriction from arbitrary $R$-linear homomorphisms to the bounded ones can be described
entirely in terms of these subobjects: $\phi$ is controlled as above precisely when there is a number $b \ge 0$ so that $\phi F(S) \subset F(S[b])$ for all choices of $S$.

In the rest of this section and the rest of the paper, we will restrict to the case of a Noetherian ring $R$.  

Let $\mathcal{P}(X)$ denote the power set of $X$ partially ordered by inclusion and viewed as a category. 
If $F$ is a left $R$-module, let $\mathcal{I}(F)$ denote the family of all $R$-submodules of $F$ partially ordered by inclusion.

\begin{DefRef}{RealBCprelim}
An $X$\textit{-filtered} $R$-module is a module $F$ together with a functor
$\mathcal{P}(X) \to \mathcal{I}(F)$
from the power set of $X$ to the family of $R$-submodules of $F$, both ordered by inclusion, such that
the value on $X$ is $F$.
It will be most convenient to think of $F$ as the functor above and use notation $F(S)$ for the value of the functor on $S$.
We will call $F$ \textit{reduced} if $F(\emptyset)=0$.

An $R$-homomorphism $f \colon F \to G$ of $X$-filtered modules is \textit{boundedly controlled} if there is a fixed number $b \ge 0$ such that the
image $f (F (S))$ is a submodule of $G (S [b])$
for all subsets $S$ of $X$.

The objects of the category $\U (X, R)$ are the reduced $X$-filtered $R$-modules, and the morphisms are the boundedly controlled homomorphisms.
\end{DefRef}

The category $\U (X, R)$ we constructed is clearly an additive category, but the more interesting structure for developing its $K$-theory is a certain Quillen exact structure.
For a good modern exposition of exact categories we refer to Keller \cite{bK:90}; there is also a leisurely review of relevant basic theory in \cite[section 2]{gCbG:00}.

Let us recall some standard terms.
If a category has kernels and cokernels for all morphisms, it is called \textit{preabelian}.  If, in addition, the canonical map $\coim (f) \to \im (f)$ for each morphism $f$ is monic and
epic but not necessarily invertible, we will say the category is
\textit{semi-abelian}, cf. \cite[pages 167-168]{wR:01} and \cite{dSsW:11}.  It is \textit{abelian} if it is also \textit{balanced} in the sense that the canonical map is an isomorphism.  
Recall also that a category is called \textit{cocomplete} if it contains colimits of
arbitrary small diagrams, cf.~Mac Lane \cite[chapter V]{sM:71}.

\begin{RemRef}{NoAb}
If $X$ is unbounded, $\U (X,R)$ is not a balanced category and therefore not an abelian category.
For an explicit description of a boundedly controlled morphism in
$\U (\mathbb{Z}, R)$ which is an isomorphism of left
$R$-modules but whose inverse is not boundedly controlled, we refer
to \cite[Example 1.5]{ePcW:85}.
\end{RemRef}

It turns out that the kernels and cokernels in $\U (X,R)$ can be characterized using an additional property a boundedly controlled
morphism may or may not have.

\begin{DefRef}{BBMor}
A morphism $f \colon F \to G$ in $\U (X,R)$ is called \textit{boundedly bicontrolled} if there exists a number $b \ge 0$ such that in addition to inclusions of submodules
\[
f (F (S)) \subset G (S [b]),
\]
there are inclusions
\[
f (F)
\cap G (S) \subset  f F (S[b])
\]
for all subsets $S \subset X$.
In this case we will say that $f$ has
\textit{filtration degree} $b$ and write $\fil (f) \le b$.
\end{DefRef}

\begin{DefRef}{GTF}
We define the
\textit{admissible monomorphisms} in $\U (X,R)$ be the boundedly
bicontrolled homomorphisms $m \colon F_1 \to F_2$ such that the map $F_1 (X) \to F_2 (X)$ is a monomorphism.
 We define the \textit{admissible epimorphisms} be the boundedly
bicontrolled homomorphisms $e \colon F_1 \to F_2$ such that $F_1 (X) \to F_2 (X)$ is an epimorphism.

Let the class $\mathcal{E}$ of
\textit{exact sequences} consist of the sequences
\[
F^{\subdot} \colon \quad F' \xrightarrow{ \ i \ } F \xrightarrow{ \ j \ } F'',
\]
where $i$ is an admissible monomorphism, $j$ is an admissible epimorphism, and $\im (i) = \ker (j)$.
\end{DefRef}

There are numerous examples of semi-abelian categories that are classical or have appeared recently in analysis and algebra that are listed in section 4 of \cite{tB:10} or in \cite{wR:01}.  For us $\U (X,R)$ is of major interest.

\begin{ThmRef}{UBPSAB}
$\U (X,R)$ is
a cocomplete semi-abelian category.
The class of exact sequences $\mathcal{E}$ gives an exact structure on $\U (X,R)$.    
\end{ThmRef}

\begin{proof}
This fact is contained in Proposition 2.6 and Theorem 2.13 of \cite{gCbG:00}. 

We want to recall the explicit construction of kernels and cokernels $\U (X,R)$ for future reference.
For a boundedly controlled morphism
$f \colon F \to G$, the kernel of $f$ in $\Mod (R)$ has the $X$-filtration $K$ where $K(S) = \ker (f) \cap F(S)$.
This gives a kernel of $f$ in $\U (X,R)$.
	 Similarly, let $I$ be the $X$-filtration of the image of
$f$ in $\Mod (R)$ by
$I(S) = \im (f) \cap G(S)$.  Let $H(X)$ be the cokernel of $f$ in $\Mod (R)$, which is the quotient $G(X)/fF(X)$.
Then there is a
filtration $H$ of $H(X)$ defined by $H(S) = \im \{ G(S)/I(S) \to H(X) \}$ where the maps between quotients are induced from the structure maps of $G$.
The resulting epimorphism $\pi \colon G(X) \to H(X)$
gives a boundedly bicontrolled morphism of filtration $0$.  Lemma 2.5 of \cite{gCbG:00} verifies the universal properties of the cokernel $H$.
\end{proof}

Now suppose $\A$ is an arbitrary cocomplete semi-abelian category.  All of the definitions and proofs presented so far can be interpreted verbatim to describe a bounded category $\U (X,\A)$ if instead of $R$-modules one uses objects from $\A$.  We should point out that the exposition in \cite{gCbG:00} is constructed in lesser generality with $\A$ assumed to be abelian.  It is important for the fibred version to relax abelian to semi-abelian. However, \cite{gCbG:00} can still serve as a good reference because throughout that paper only semi-abelian properties of $\A$ are used.  

In particular we have this conclusion.

\begin{ThmRef}{UBPSAB2}
For any cocomplete semi-abelian category $\A$, $\U (X,\A)$ is
a cocomplete semi-abelian category.   
\end{ThmRef}

Of course, the category of $R$-modules $\Mod (R)$ is a cocomplete abelian category and so can serve as a basic example of cocomplete semi-abelian coefficients $\A$.  In this case $\U (X,\A)$ is precisely $\U (X,R)$.  

\SSecRef{Properties of filtered objects}{PFO}

In this section we will start imposing several conditions on objects in the bounded category $\U (X,\A)$.  These conditions are agnostic to the nature of $\A$, and so we can use the simpler notation $\U (X)$ for this category.  These conditions and results about them will be used in the contexts of both module categories and abstract semi-abelian categories as coefficients.

\begin{DefRef}{HYUT}
Let $F$ be an $X$-filtered $R$-module.
\begin{itemize}
\item $F$ is called \textit{lean} or $D$-\textit{lean} if there is a number $D \ge 0$ such that
\[
F(S) \subset \sum_{x \in S} F(x[D])
\]
for every subset $S$ of $X$.
\item $F$ is called \textit{split} or $D'$-\textit{split} if there is a
number $D' \ge 0$ such that we have
\[
F(S) \subset F(T[D']) + F(U[D'])
\]
whenever a subset $S$ of $X$ is written as a union $T \cup U$.
\item $F$ is called \textit{insular} or $d$-\textit{insular} if there is a
number $d \ge 0$ such that
\[
F(S) \cap F(U) \subset F(S[d] \cap U[d])
\]
for every pair of subsets $S$, $U$ of $X$.
\end{itemize}
\end{DefRef}

\begin{PropRef}{leaninscld}
The properties of being lean, split, and insular are preserved under
isomorphisms in $\U (X)$.  Also, a $D$-lean filtered module is $D$-split. 
\end{PropRef}

\begin{proof}
If $f \colon F_1 \to F_2$ is an isomorphism with $\fil (f) \le b$, and $F_1$ is $D$-lean, $D'$-split, and $d$-insular, then $F_2$ is $(D+b)$-lean, $(D'+b)$-split, and $(d+2b)$-insular.  
For the other statement, we have
\[
F(T \cup U) \subset \sum_{x \in T} F(x[D]) + \sum_{x \in U} F(x[D]) \subset F(T[D]) + F(U[D])
\]
since in general $\sum_{x \in S} F(x[D]) \subset F(S[D])$.
\end{proof}

A collection of objects in an exact category is said to be \textit{closed under extensions} if the middle term of an exact sequence belongs to the collection in case both of the extreme terms belong to the collection.

\begin{LemRef}{LandS}
$\mathrm{(1)}$ Lean objects are closed under extensions.

$\mathrm{(2)}$ Insular objects are closed under extensions.

$\mathrm{(3)}$ Split objects are closed under extensions.
\end{LemRef}

\begin{proof}
For an exact sequence
$E' \xrightarrow{f} E \xrightarrow{g} E''$ in $\U (X)$, let $b \ge 0$ be a common
filtration degree for $f$ and $g$.
The first two statements follow from parts (1) and (2) of \cite[Proposition 2.18]{gCbG:00}.  It is shown there that if $E'$ and $E''$ are $D$-lean then $E$ is $(4b+D)$-lean. Also, if both $E'$ and $E''$ are $d$-insular then $E$ is $(4b + 2d)$-insular.

To prove (3), suppose both $E'$ and $E''$ are $D'$-split.
We have
\[
gE(T \cup U) \subset E'' (T[b] \cup U[b]),
\]
because in general $(T \cup U)[b] \subset T[b] \cup U[b]$.
So
\begin{equation} \begin{split}
&g (T \cup U)\\
\subset\ &E'' (T[b + D']) + E'' (U[b + D'])\\
\subset\ &gE (T[2b + D']) + gE (U[2b + D']).
\end{split} \notag \end{equation}
If $z \in E(T \cup U)$ then we can write $g(z) = g(z_1) + g(z_2)$ where $z_1 \in E(T[2b + D'])$ and $z_2 \in E(U[2b + D'])$.
Since $z - z_1 - z_2$ is an element of $\ker (g) \cap E(T[2b + D'] \cup U[2b + D'])$, we have an element
\begin{equation} \begin{split}
k \in\  &E' (T[3b + D'] \cup U[3b + D'])\\
\subset\ &E' (T[3b + 2D']) + E' (U[3b + 2D'])
\end{split} \notag \end{equation}
such that
\[
z = f(k) + z_1 + z_2 \in E (T[4b + 2D']) + E (U[4b + 2D']).
\]
So $E$ is $(4b + 2D')$-split.
\end{proof}

\begin{LemRef}{lninpres}
Let
\[
E' \xrightarrow{\ f \ } E \xrightarrow{\ g \ } E''
\]
be an exact sequence in $\U (X)$.

\begin{enumerate}
\item If the object $E$ is
lean then $E''$ is lean.

\item If $E$ is
split then $E''$ is split.

\item If $E$ is
insular then $E'$ is insular.

\item If $E$ is
insular and $E'$ is lean then $E''$ is insular.

\item If $E$ is
insular and $E'$ is split then $E''$ is insular.

\item If $E$ is
split and $E''$ is insular then $E'$ is split.
\end{enumerate}
\end{LemRef}

\begin{proof}
Let $b \ge 0$ be a common
filtration degree for $f$ and $g$.
If $E$ is $D$-lean, $D'$-split, or $d$-insular, it is easy to show that $E''$ is $(D+2b)$-lean or $(D'+2b)$-split and $E'$ is $(d+2b)$-insular respectively, which verifies (1), (2), and (3).

Statement (4) follows from the proof of part (3c) of \cite[Proposition 2.18]{gCbG:00}.  It is shown there that if $E'$ is $D$-lean and $E$ is $d$-insular then $E''$ is $(4b + D + d)$-insular.  The same proof actually shows statement (5).  The only equation in that proof that uses $D$-leanness of $E'$ is only used to get a consequence that is in fact immediate from the assumption that $E'$ is $D$-split.

(6)  Suppose $E$ is $D'$-split and $E''$ is $d$-insular.
Given $z \in E' (T \cup U)$, we have $f(z) \in E (T[b] \cup U[b])$.
Now 
$f(z) \in E (T[b + D']) + E (U[b + D'])$, so we can write accordingly $f(z) = y_1 + y_2$.
Now $f(z) \in \ker (g)$, because $g(y_1) + g(y_2) = 0$.
Since $E''$ is $d$-insular,
\[
g(y_1) = - g(y_2) \in E'' (T[2b + D' + d] \cap U[2b + D' + d]),
\]
so we are able to find
\[
y \in E (T[3b + D' + d] \cap U[3b + D' + d])
\]
such that $g(y) = g(y_1) = - g(y_2)$, because generally $(S \cap P)[b] \subset S[b] \cap P[b]$.
Thus
\[
f(z) = y_1 + y_2 = (y_1 - y) + (y_2 + y)
\]
and
\[
y_1 - y \in E(T[3b + D' + d]),
\quad
y_2 + y \in E(U[3b + D' + d]).
\]
Let $z_1 = f^{-1} (y_1 - y)$ and $z_2 = f^{-1} (y_2 + y)$, and we have $z = z_1 + z_2$ such that
\[
z_1 \in E' (T[4b + D' + d]),
\quad
z_2 \in E' (U[4b + D' + d]),
\]
so $E'$ is $(4b + D' + d)$-split.
\end{proof}

\begin{CorRef}{WOW}
Let
$E' \xrightarrow{f} E \xrightarrow{g} E''$
be an exact sequence in $\U (X)$.
If $E$ is
split and insular then $E''$ is insular if and only if $E'$ is split.
\end{CorRef}

\begin{proof}
This fact is the combination of parts (5) and (6) of the Lemma.
\end{proof}

\begin{RemRef}{WOW2}
The last Corollary is in contrast with the absence of the analogous general fact if one substitutes the lean property for the split property.
However, the analogue is true in the presence of certain geometric assumptions on the metric space.
For example, suppose $X$ has finite asymptotic dimension.
Then from the main theorem of \cite{gCbG:03}, we have the following counterpart to part (6) of the Lemma:
if $E$ is lean and $E''$ is insular then $E'$ is lean.
This fact is not needed in this paper.
Here, the excision properties of the theory rely only on the properties of the cokernels.
For the applications in \cite{gCbG:15}, properties of the kernels become crucial in dealing with coherence issues, and the geometric conditions need to be imposed.
\end{RemRef}

\begin{DefRef}{TRYE}
We define $\LI (X)$ as the full subcategory of $\U (X)$ on objects that are lean and insular with the induced exact structure.
Similarly, $\SI (X)$ is the full subcategory of $\U (X)$ on objects that are split and insular.
\end{DefRef}

Exact structures in $\LI (X)$ and $\SI (X)$ can be induced from $\U (X)$.
A full subcategory $\mathbf{H}$ of an exact category $\C$ is said to
be \textit{closed under extensions} or \textit{thick} in $\C$ if 
\begin{enumerate}
	\item $\mathbf{H}$ contains the
zero object, and
\item for any exact sequence $C' \to C \to C''$ in $\C$,
if $C'$ and $C''$ are isomorphic to objects from $\mathbf{H}$ then so is
$C$.
\end{enumerate}
It is known (cf. \cite[Lemma 10.20]{tB:10}) that a subcategory closed under
extensions in $\C$ inherits the exact structure from $\C$.

\begin{ThmRef}{OneTwo}
$\LI (X)$ and $\SI (X)$ are closed under extensions in $\U (X)$.
Therefore, $\LI (X)$ and $\SI (X)$ are exact subcategories of $\U (X)$,
so we have a sequence of exact inclusions
\[
\LI (X) \longrightarrow \SI (X) \longrightarrow \U (X).
\]
\end{ThmRef}

\begin{proof}
The first fact follows from parts (1) and (2) of Lemma \refL{LandS}, the second from (2) and (3).
\end{proof}

\SSecRef{Local finiteness property}{LF}

Finally, there is an additional property that will consider only in module categories.

\begin{DefRef}{RealBCprel}
An $X$-filtered $R$-module $F$ is \textit{locally finitely generated} if $F (S)$ is a finitely generated $R$-module for
every bounded subset $S \subset X$.

The category $\BLI (X,R)$ is the full
subcategory of $\LI (X,R)$ on the locally finitely generated objects.
Similarly, the companion category $\BSI (X,R)$ is the full
subcategory of $\SI (X,R)$ on the locally finitely generated objects.
\end{DefRef}

\begin{ThmRef}{ExtCl}
The category $\BLI (X,R)$ is closed under extensions in $\LI (X,R)$.
Similarly, the category $\BSI (X,R)$ is closed under extensions in $\SI (X,R)$.
\end{ThmRef}

\begin{proof}
If $f \colon F \to G$ is an isomorphism with $\fil (f) \le b$ and $G$ is locally finitely generated, then $F (U)$ are finitely generated submodules of $G (U[b])$ for all bounded $U$, since $R$ is a Noetherian ring.
Suppose
\[
F' \xrightarrow{\ f \ } F \xrightarrow{\ g \ } F''
\]
is an exact sequence and let $b \ge 0$ be a common
filtration degree for both $f$ and $g$.
Assume that $F'$ and $F''$ are locally finitely generated.
For every bounded
subset $U \subset X$ the restriction $g \colon F(U) \to
gF(U)$ is an epimorphism onto a submodule of the finitely generated $R$-module
$F''(U[b])$. The kernel of $g \vert F(U)$ is a submodule of $F'(U[b])$, which is also
finitely generated.  So the extension $F(U)$ is finitely generated.
\end{proof}

\begin{CorRef}{BisWAb}
$\BLI (X,R)$ and $\BSI (X,R)$ are exact categories.  The
additive category $\mathcal{C}(X,R)$ of geometric $R$-modules with
the split exact structure is an exact subcategory of $\BLI (X,R)$, so there is a sequence of exact inclusions
\[
\mathcal{C}(X,R) \longrightarrow \BLI (X,R) \longrightarrow \BSI (X,R) \longrightarrow \U (X,R).
\]
\end{CorRef}

\begin{RemRef}{REDCVB}
We want to briefly explain the roles played by the two conditions, lean and split, that distinguish the two categories $\BLI (X,R)$ and $\BSI (X,R)$.  $\BLI (X,R)$ was used exclusively in \cite{gCbG:00}, where it was proven to have good excision properties.  There is a separate important issue of homological coherence that still requires the lean condition for its resolution, cf. \cite{gCbG:03,gCbG:15}.
The setting with the split condition in $\BSI (X,R)$ is much more streamlined for the excision arguments but has insufficient coherence properties.  
In the next section, we will pursue the goal of combining the two different conditions in the ``base'' and ``fibre'' in order to achieve required coherence in the base and ``fibrewise'' excision properties.  The hybrid lean/split condition will provide a considerable advantage because the fibred setting is more complicated than the absolute case.  
\end{RemRef}

Recall that a morphism $e \colon F \to F$ is an idempotent if
$e^2=e$. Categories in which every idempotent is the projection
onto a direct summand of $F$ are called \textit{idempotent complete}.

\begin{PropRef}{PAIC}
$\BLI (X,R)$ and $\BSI (X,R)$ are idempotent
complete.
\end{PropRef}

\begin{proof}
First note that a regular preabelian category is idempotent complete.  The proof is exactly the same as for an abelian category: if $e$ is an idempotent then its kernel is split by $1-e$.
Since the restriction of an idempotent $e$ to the image of $e$ is
the identity, every idempotent here is boundedly bicontrolled of
filtration $0$. It follows easily that the splitting of $e$ in
$\Mod(R)$ is in fact a splitting in $\BLI (X,R)$ or $\BSI (X,R)$.
\end{proof}

Finally, we need to address (the lack of) inheritance in filtered modules.  It is immediate that a submodule of an insular filtered module is also insular with respect to the standard filtration induced on the submodule.  However, a submodule of a lean filtered module is not necessarily lean.

\begin{DefRef}{Strictness}
An $X$-filtered object $F$ is called
\textit{strict} if there exists an order
preserving function
$\ell \colon \mathcal{P}(X) \to [0, +\infty)$
such that for every $S \subset X$ the submodule $F(S)$ is
$\ell_S$-lean and $\ell_S$-insular with respect to the standard
$X$-filtration $F(S)(T) = F(S) \cap F(T)$.
\end{DefRef}

It is important to note that this property is not preserved under isomorphisms, so the subcategory of strict
objects is not essentially full in $\BLI (X,R)$.

The bounded category $\B (X,R)$ was defined in \cite{gCbG:00} as the full subcategory of $\BLI (X,R)$ on objects isomorphic to strict objects. Now this category is closed under exact extensions in $\U (X,R)$ according to \cite[Theorem 2.22]{gCbG:00} and so is an exact category. 

A consequence of strictness, or more generally being isomorphic to a strict object, is the following feature.
Given a filtered module $F$ in $\B (X, R)$, a \textit{lean grading} of $F$ is a functor
$\widetilde{F} \colon \mathcal{P}(X) \to \mathcal{I}(F)$
from the power set of $X$ to the submodules of $F$ such that
\begin{enumerate}
\item each $\widetilde{F} (S)$ is an object of $\BLI (X, R)$ when given the standard filtration,
\item there is a number $K \ge 0$ such that
\[
F(S) \subset \widetilde{F} (S) \subset F(S[K])
\]
for all subsets $S$ of $X$.
\end{enumerate}
Clearly, each $\widetilde{F} (S)$ is an object of $\B (X,R)$. 
Also an actual strict object has a lean grading by $\widetilde{F} (S) = F(S)$ with $K = 0$.

We note for the interested reader that the theory in \cite{gCbG:00}, including the excision theorems, could be alternately developed for modules with lean gradings in place of $\B (X,R)$.
We do not require such theory in this paper.
Instead, we develop a similar but more relaxed notion of gradings in $\BSI (X,R)$.

\SSecRef{Graded objects and their closure properties}{GrCP}

\begin{DefRef}{WEI}
Given a filtered module $F$ in $\BSI (X, R)$, a \textit{grading} of $F$ is a functor
$\mathcal{F} \colon \mathcal{P}(X) \to \mathcal{I}(F)$ such that
\begin{enumerate}
\item each $\mathcal{F} (S)$ is an object of $\BSI (X, R)$ when given the standard filtration,
\item there is a number $K \ge 0$ so that
\[
F(S) \subset \mathcal{F} (S) \subset F(S[K])
\]
for all subsets $S$ of $X$.
\end{enumerate}
We will say that a filtered module $F$ is \textit{graded} if it is possible to equip it with a grading, but there is no specific choice of grading that is specified.
\end{DefRef}

\begin{PropRef}{OINBFG}
The graded objects are closed under isomorphisms.
\end{PropRef}

\begin{proof}
If $f \colon F \to F'$ is an isomorphism and $F$ has a grading $\mathcal{F}$, a grading for $F'$ is given by
$\mathcal{F}' (C) \, = \, f \mathcal{F} (C[K+b])$,
where $b$ is a filtration bound for $f$.
\end{proof}

\begin{DefRef}{JHSD}
We define $\G (X,R)$ as the full subcategory of $\BSI (X,R)$ on the locally finitely generated graded filtered modules.
\end{DefRef}

\begin{PropRef}{VFGHJ}
$\G (X,R)$ is closed under extensions in $\BSI (X,R)$.
Therefore $\G (X,R)$ is an exact subcategory of $\BSI (X,R)$.
\end{PropRef}

\begin{proof}
Given an exact sequence $F \xrightarrow{f} G \xrightarrow{g} H$ in $\BSI (X,R)$,
let $b \ge 0$ be a common
filtration degree for both $f$ and $g$ as boundedly bicontrolled
maps, and assume that
$F$ and $H$ are graded modules in $\G (X,R)$ with the associated functors $\mathcal{F}$ and $\mathcal{H}$.

To define a grading for $G$,
consider a subset $S$ and suppose $\mathcal{H} (S[b])$
is $D$-split and $d$-insular. The induced epimorphism
$g \colon G(S[2b]) \cap g^{-1} \mathcal{H} (S[b]) \to 
\mathcal{H} (S[b])$
extends to another epimorphism
\[
g' \colon f \mathcal{F} (S[3b]) + G(S[2b]) \cap g^{-1} \mathcal{H} (S[b])
\longrightarrow \mathcal{H} (S[b])
\]
with $\ker (g') = \mathcal{F} (S[3b])$.
Without loss of generality, suppose $\mathcal{F} (S[3b])$ is $D$-split and $d$-insular.
We define
\[
\mathcal{G}(S) = f \mathcal{F} (S[3b]) + G(S[2b]) \cap g^{-1} \mathcal{H} (S[b]).
\]
From parts (2) and (3) of Lemma \refL{LandS}, the module $\mathcal{G}(S)$ with the standard filtration is $(4b+2d)$-split and $(4b+2d)$-insular.
Since $G(S) \subset g^{-1} \mathcal{H} (S[b])$, we have $G(S) \subset \mathcal{G}(S)$.
On the other hand, if the grading $\mathcal{F}$ has characteristic number $K \ge 0$ then
$\mathcal{G}(S) \subset G(S[4b+K])$.
The last fact together with Theorem \refT{ExtCl} shows that $\mathcal{G}(S)$ is finitely generated.
\end{proof}

The relations between the categories in this section can be summarized as a commutative diagram of fully exact inclusions
\begin{equation}
\xymatrix@R=3mm@C=8mm{
 &\BLI (X,R) \ar[r]
 &\BSI (X,R) \\
  \mathcal{C} (X,R) \ar@/^/[ur] \ar@/_/[dr] \\
 &\B (X,R) \ar[r] \ar[uu]
 &\G (X,R) \ar[uu]
} \notag
\end{equation}

The advantage of working with the category $\G (X,R)$ is that one can readily localize to geometrically defined subobjects.

\begin{LemRef}{JHQASF}
Suppose $G$ is a graded $X$-filtered module with a grading $\mathcal{G}$.  Let $F$ be a submodule which is split with respect to the standard filtration.  Then $\mathcal{F} (S) = F \cap \mathcal{G} (S)$ is a grading of $F$.
\end{LemRef}

We will call the grading of a split submodule $F$ obtained in Lemma \refL{JHQASF} the \textit{standard grading} of the submodule.

\begin{proof}
Of course, $F(S) = F \cap G(S) \subset F \cap \mathcal{G} (S) = \mathcal{F}(S)$.
On the other hand, there is $d \ge 0$ such that $\mathcal{G} (S) \subset G(S[d])$, so $\mathcal{F} (S) \subset F \cap G(S[d]) = F(S[d])$.

Consider the inclusion of modules $i \colon F \to G$, and take the quotient $q \colon G \to H$.
Both $F$ and $G$ are split and insular, so $H$ is split and insular by parts (2) and (4) of Lemma \refL{lninpres}, with respect to the quotient filtration.
We define $\mathcal{H}(S)$ as the partial image $q \mathcal{G}(S)$ and give $\mathcal{H}(S)$ the standard filtration in $H$.  Then $\mathcal{H}(S)$ is split as the image of a split $\mathcal{G}(S)$ and insular since $H$ is insular.  Now the kernel of the epimorphism $q \vert \colon \mathcal{G}(S) \to \mathcal{H}(S)$,
which is $F \cap \mathcal{G}(S)$, is split by part (6) of Lemma \refL{lninpres}.
Since $F$ is insular, $\mathcal{F}(S)$ is also insular.
This shows that $\mathcal{F}(S)$ gives a grading for $F$.
\end{proof}
 
 This can be promoted to the following result.

\begin{PropRef}{HUVZFMN}
Given a boundedly bicontrolled epimorphism $g \colon G \to H$ in $\BSI (X,R)$, suppose $F$ is a submodule of $G$ which is the kernel of $g$ in $\Mod (R)$.  It is given the standard filtration. If $G$ is graded and $F$ is split then both $H$ and $F$ are graded.
\end{PropRef}

\begin{proof}
The grading for $H$ is given by $\mathcal{H}(S) = g \mathcal{G}(S[b])$, where $b$ is a chosen bicontrol bound for $g$.
Each $\mathcal{H}(S)$ is split and insular as in the proof of Lemma \refL{JHQASF}.
The inclusions $H(S) \subset gG(S[b]) \subset g \mathcal{G}(S[b]) = \mathcal{H}(S)$ and
$g \mathcal{G}(S[b]) \subset g {G}(S[b + K]) \subset H(S[2b+K])$ show that $\mathcal{H}$ is a grading.
The same argument as in Lemma \refL{JHQASF} shows that $\mathcal{F}(S) = F \cap \mathcal{G}(S[b])$ gives a grading for $F$.
\end{proof}

Applying this fact, we are able to characterize admissible monomorphisms in $\G(X,R)$ as follows.

\begin{PropRef}{NJWSX}
The inclusion of a subobject $i \colon F \to G$ in $\G(X,R)$ is an admissible monomorphism if and only if $F$ is split.	
\end{PropRef}

\begin{proof}
The cokernel $H$ of $i$ 	in $\U (X,R)$ has the filtration described in the proof of Theorem \refT{UBPSAB}. From parts (2) and (5) of Lemma \refL{lninpres}, $H$ is split and insular.  In fact, $H$ is a cokernel of $i$ in $\BSI (X,R)$.
From Proposition \refP{HUVZFMN}, $H$ is graded, so it is also a cokernel of $i$ in $\G(X,R)$.
\end{proof}

We will use the following convention.
When $d \le 0$, the notation $S[d]$ will stand for the subset $S \setminus (X \setminus S)[-d]$.

\begin{CorRef}{LPCSL}
Given an object $F$ in $\G (X,R)$ and a subset $S$ of $X$, there is a number $K \ge 0$ and an admissible subobject $i \colon F_S \to F$ in $\G (X,R)$ with the property that $F_S \subset F(S[K])$.
Moreover, the cokernel $q \colon F \to H$ has the property that $H (X) = H ((X \setminus S) [2D'])$, where $D'$ is a splitting constant for $F$.
\end{CorRef}

\begin{proof}
The object $F$ has a grading $\mathcal{F}$ with a characteristic constant $K$.
Property (1) in Definition \refD{WEI} guarantees that $\mathcal{F}(S)$ is a split object for any subset $S$.
For the first statement, choose $F_S = \mathcal{F}(S)$ with the grading defined in Lemma \refL{JHQASF} and apply Proposition \refP{NJWSX}. 
The second statement is shown as follows.
By part (2) of Lemma \refL{lninpres}, since $\fil (q) = 0$, if $F$ is $D'$-split then $H$ is $D'$-split.
Let $T = S[-D']$, then $T[D'] \subset S$, so 
\[
H(T[D']) = qF(T[D']) \subset qF(S) \subset qF_S = 0.
\]
Using the decomposition $X = T \cup (X \setminus T)$ we can write
\[
H(X) = H(T[D']) + H((X \setminus T)[D']) = H((X \setminus T)[D']) = H((X \setminus S)[2D']).  \qedhere
\] 
\end{proof}

The last three results can be summarized as follows.

\begin{ThmRef}{LPCSLii}
Given a graded object $F$ in $\G (X,R)$ and a subset $S$ of $X$,
we assume that $F$ is $D'$-split and $d$-insular and is graded by $\mathcal{F}$.
The submodules $\mathcal{F} (S)$ have the following properties:
\begin{enumerate}
\item each $\mathcal{F} (S)$ is graded by $\mathcal{F}_S (T) = \mathcal{F} (S) \cap \mathcal{F} (T)$;
\item $F(S) \subset \mathcal{F} (S) \subset F(S[K])$ for some fixed number $K \ge 0$;
\item suppose $q \colon F \to H$ is the cokernel of the inclusion $i \colon \mathcal{F} (S) \to F$, then $H$ is supported on $(X \setminus S) [2D']$;
\item $H (S[-2D' -2d]) =0$.
\end{enumerate}
\end{ThmRef}

\begin{proof}
Properties (1), (2), (3) are consequences of the last four results.
(4) follows from the fact that a $d$-insular filtered module is $2d$-separated,
in the sense that for any pair of subsets $S$ and $T$ such that $S[2d] \cap T = \emptyset$ we have
$S[d] \cap T[d] = \emptyset$ so $F(S) \cap F(T) = 0$.
Now $H(S[-2D' -2d]) \cap H((X \setminus S)[2D']) =0$, but
$H((X \setminus S)[2D']) = H(X)$, thus $H(S[-2D' -2d]) = 0$.
\end{proof}

\begin{RemRef}{NOFUNCT}
	Functoriality properties in controlled theories are well-understood.  As expected, bounded $G$-theory is covariantly functorial in both variables in the sense that $\G (X,R)$ is a covariant functor from the category  of proper metric spaces and uniformly expansive maps to exact categories when $R$ is fixed and is similarly a covariant functor from Noetherian rings to exact categories when $X$ is fixed.  These details become important in the construction of the equivariant theory and in specific applications.  We avoid questions of functoriality in this paper as we concentrate on computational tools such as excision.  
\end{RemRef}

\SecRef{Fibred bounded \textit{G}-theory}{MSOGS2}

\SSecRef{Introduction of fibred control in \textit{G}-theory}{IFC}

Suppose $X$ and $Y$ are two proper metric spaces and $R$ is a Noetherian ring. 
The product $X \times Y$ is given the product metric
$d((x,y), (x',y')) = \max \{ d(x,x'),d(y,y') \}$.
There is certainly the semi-abelian category $\U (X \times Y, R)$, the exact category $\LI (X \times Y, R)$ and, further, the bounded category $\BLI (X \times Y, R)$.

We wish to construct a larger \textit{fibred bounded category} $\BX (Y)$.
The result will involve a mix of features from $\BLI (X, \A)$ and $\BSI (Y,R)$ and contain $\BLI (X \times Y, R)$ as an exact subcategory.

One definition can be made by simply imitating $K$-theory with fibred control as described in Remark \refR{HAVRCF}.  It is obtained as an iterated construction from the end of section \refSS{EBKT} on the level of \textit{unrestricted} category as $\U (X, \U(Y,R))$.
From Theorem \refT{UBPSAB2}, $\U (X, \U(Y,R))$ is a complete semi-abelian category with complete semi-abelian coefficients $\A = \U(Y,R)$.  This concept has a transparent definition but is not well-suited for a fibred theory.  In particular, while it does contain $\mathcal{C} (X \times Y, R)$, it does not contain the fibred bounded category $\mathcal{C} (X, \mathcal{C}(Y,R))$ as a subcategory.  Now we proceed to develop a different, more explicit set-up in terms of a category $\U_X (Y)$ that naturally contains $\U (X \times Y, R)$ and $\mathcal{C} (X, \mathcal{C}(Y,R))$ as exact subcategories.  

\begin{DefRef}{GEPCW2}
Given an $R$-module $F$, an $(X,Y)$\textit{-filtration} of $F$ is a functor
$ \phi_F \colon \mathcal{P}(X \times Y) \to  \mathcal{I}(F)$
from the power set of the product metric space to the partially ordered family of $R$-submodules of $F(X \times Y)$.
Whenever $F$ is given a filtration, and there is no ambiguity, we will denote the values $\phi_F (U)$ by $F(U)$.
We assume that $F$ is \textit{reduced}
in the sense that the value on the empty subset is $0$.

The associated $X$-filtered $R$-module $F_X$ is given by \[ F_X (S) = F(S \times Y).\]  Similarly, for each subset $S \subset X$, one has the $Y$-filtered $R$-module $F^S$ given by \[ F^S (T) = F(S \times T). \]
In particular, $F^X(T) = F(X \times T)$.
\end{DefRef}

We will use the following notation generalizing enlargements in a metric space.

\begin{NotRef}{MIONBZ}
Given a subset $U$ of $X \times Y$ and a function $k \colon X \to [0, + \infty )$, let
\[
U [k] = \{ (x,y) \in X \times Y \ \vert \ \textrm{there\ is} \ (x,y') \in U \ \textrm{with} \ d(y,y') \le k(x) \}.
\]
If in addition we are given a number $K \ge 0$ then
\[
U [K,k] = \{ (x,y) \in X \times Y \ \vert \ \textrm{there\ is} \ (x',y) \in U[k] \ \textrm{with} \ d(x,x') \le K  \}.
\]
So $U[k] = U[0,k]$.  Notice that
if $U$ is a single point $(x,y)$ then
\[
U[K,k] = x[K] \times y[k(x)] = (x,y)[K,0] \times (x,y)[0,k(x)].
\]
More generally, one can equivalently write
\[
U[K,k] = \bigcup_{(x,y) \in U} x[K] \times y[k(x)].
\]
If $U$ is a product set $S \times T$, it will be convenient to use the notation $(S,T)[K,k]$ in place of $(S \times T)[K,k]$.
More generally, because the roles of the factors are very different when working with $(X,Y)$-filtrations, we will use the notation $(X,Y)$ for the product metric space so that the order of the factors is unambiguous.  Similarly, we will use the notation $(S,T)$ for the product subset $S \times T$ in $(X,Y)$.
\end{NotRef}

\begin{DefRef}{LOIV}
We will refer to the pair $(K,k)$ in the notation $U[K,k]$ as the \textit{enlargement data}.
\end{DefRef}

It is clear that when $Y = \point$, $U [K,k] = U [K]$ for any function $k$ under the identification $X \times Y = X$.

\begin{NotRef}{UYBV}
Let $x_0$ be a chosen fixed point in $X$.
Given a monotone function $h \colon [0, + \infty ) \to [0, + \infty )$, there is a function $h_{x_0} \colon X \to [0, + \infty )$ defined by
\[
h_{x_0} (x) = h (d_X (x_0,x)).
\]
\end{NotRef}

\begin{DefRef}{EqBC}
Given two $(X,Y)$-filtered modules $F$ and $G$, an $R$-homomorphism $f \colon F(X \times Y) \to G(X \times Y)$ is \textit{boundedly controlled} if
there are a number $b \ge 0$ and a monotone function $\theta \colon [0, + \infty ) \to [0, + \infty )$ such that
\begin{equation}
fF(U) \subset G(U[b,\theta_{x_0}]) \tag{$\dagger$}
\end{equation}
for all subsets $U \subset X \times Y$ and some choice of $x_0 \in X$.
It is clear that the condition is independent of the choice of $x_0$.

The \textit{unrestricted fibred bounded category} $\U_X (Y)$ has $(X,Y)$-filtered modules as objects and the boundedly controlled homomorphisms as morphisms.
\end{DefRef}

\begin{ThmRef}{MNBVX}
$\U_X (Y)$ is a cocomplete semi-abelian category.	
\end{ThmRef}

First we require a very useful fact.

\smallskip

A morphism $f \colon F \to G$ in $\U_X (Y)$ is \textit{boundedly bicontrolled} if
there is filtration data $b \le 0$ and $\theta \colon [0, + \infty ) \to [0, + \infty )$ as in Definition \refD{EqBC}, and
in addition to ($\dagger$) one also has the containments
$fF \cap G(U) \subset fF(U[b,\theta_{x_0}])$.
In this case, we will use the notation $\fil (f) \le (b,\theta)$.

\begin{LemRef}{CatOK}
Let $f_1 \colon F \to G$, $f_2 \colon G \to H$ be in $\U_X (Y)$
and $f_3 = f_2  f_1$.
\begin{enumerate}
\item If $f_1$, $f_2$ are boundedly bicontrolled morphisms
and either $f_1 \colon F(X \times Y) \to G(X \times Y)$ is an epi or $f_2 \colon
G(X \times Y) \to H(X \times Y)$ is a monic, then $f_3$ is also boundedly
bicontrolled.
\item If $f_1$, $f_3$ are boundedly bicontrolled
and $f_1$ is epic then $f_2$ is also boundedly bicontrolled; if
$f_3$ is only boundedly controlled then $f_2$ is also boundedly
controlled.
\item If $f_2$, $f_3$ are boundedly bicontrolled
and $f_2$ is monic then $f_1$ is also boundedly bicontrolled; if
$f_3$ is only boundedly controlled then $f_1$ is also boundedly
controlled.
\end{enumerate}
\end{LemRef}

\begin{proof}
Suppose $\fil (f_i) \le (b,\theta)$ and $\fil (f_j) \le (b',\theta')$ for $\{ i,j \}
\subset \{ 1, 2, 3 \}$, then in fact $\fil (f_{6-i-j}) \le (b+b',\theta+\theta')$ in each of the three cases. For example, there are
factorizations
\begin{gather}
f_2 G(U) \subset f_2 f_1 F(U[b,\theta_{x_0}]) =
f_3 F(U[b,\theta_{x_0}]) \subset H(U[b+b', \theta_{x_0} + \theta'_{x_0}]) \notag \\
f_2 G(X) \cap H(U) \subset f_3 F(U[b',\theta'_{x_0}]) = f_2 f_1 F(U[b',\theta'_{x_0}])
\subset f_2 G(U[b+b', \theta_{x_0} + \theta'_{x_0}]) \notag
\end{gather}
which verify part 2 with $i=1$, $j=3$.
\end{proof}

\begin{proof}[Proof of Theorem \refT{MNBVX}]
The additive properties are inherited from $\Mod (R)$, so the biproduct is given by the filtration-wise operation
\(
(F \oplus G)(U)
= F(U) \oplus G(U)
\)
in $\Mod (R)$. For any boundedly controlled morphism
$f \colon F \to G$, the kernel of $f$ in $\Mod (R)$ has the standard
$(X,Y)$-filtration $K$ where
\(
K(S) = \ker (f) \cap F(S)
\)
which gives the
kernel of $f$ in $\U_X (Y)$. The canonical monic $\kappa \colon K
\to F$ has filtration data $(0,0)$ and is therefore boundedly bicontrolled.
It follows from part 3 of Lemma
\refL{CatOK} that $K$ has the universal properties of the kernel in
$\U_X (Y)$.

Similarly, let $I$ be the standard $(X,Y)$-filtration of the image of
$f$ in $\Mod(R)$ by
\(
I(U) = \im (f) \cap G(U).
\)
Then there is a
presheaf $C$ over $(X,Y)$ with
\(
C(U) = G(U)/I(U)
\)
for $U \subset (X,Y)$.
Of course $C(X \times Y)$ is the cokernel of $f$ in $\Mod (R)$. Consider an $(X,Y)$-filtered object $\overline{C}$ associated to $C$ given by
\(
\overline{C} (U)
= \im C (U, X \times Y).
\)
The canonical morphism $\pi \colon G(X \times Y) \to C(X \times Y)$
gives a boundedly bicontrolled morphism $\pi \colon G \to \overline{C}$ of filtration $(0,0)$  since
\[
\im (\pi  G(U,X \times Y)) = \im C(U,X \times Y) = \overline{C} (U).
\]
This in conjunction with part 2 of Lemma \refL{CatOK} also verifies
the universal cokernel properties of $\overline{C}$ and $\pi$ in $\U_X
(Y)$.
\end{proof}

We mention one useful perspective on $\U_X (Y)$.  

\begin{PropRef}{HUNBV3}
Suppose $f \colon F \to G$ in $\U_X (Y)$ is boundedly controlled with control data $(b, \theta)$.
Then
\begin{enumerate}
\item $f$ is bounded by $b$ when viewed as a morphism $F_X \to G_X$ in $\U (X, R)$, and
\item for each bounded subset $S \subset X$, the restriction
$f \vert \colon F_X (S) \to G_X (S[b])$
is bounded when viewed as a morphism $F^S \to G^{S[b]}$ of $Y$-filtered modules in $\U (Y)$.
\end{enumerate}
\end{PropRef}

\begin{proof}
If $f \colon F \to G$ is $(b,\theta)$-controlled then for any subset $S \subset X$ we have $f F_X (S) \subset
G((S,Y)[b,\theta_{x_0}]) \subset G(S[b],Y) = G_X (S[b])$.  So $f \colon F_X \to G_X$ is bounded by $b$.
Now for a given subset $S \subset X$, let us define $\theta_S = \sup_{x \in S} \theta_{x_0} (x)$.
Then $fF_X (S)(T) = fF(S,T) \subset G(S[b],T[\theta_S]) = G_X (S[b])(T[\theta_S])$
verifying that $f \vert \colon F^S \to G^{S[b]}$ is bounded by $\theta_S$.
\end{proof}

\SSecRef{Properties of fibred objects}{HIAHG}

\begin{DefRef}{StrFib}
An $(X,Y)$-filtered module $F$ is called
\begin{itemize}
\item \textit{lean} or $(D,\Delta)$-\textit{lean} if there is a number $D \ge 0$ and a monotone function
$\Delta \colon [0, + \infty ) \to [0, + \infty )$ so that
    \begin{align}
    F(U) &\subset \sum_{(x,y) \in U} F(x[D] \times y[\Delta_{x_0}(x)]) = \sum_{(x,y) \in U} F((x,y)[D,\Delta_{x_0}])  \notag
    \end{align}
    for any subset $U$ of $X \times Y$, \medskip
\item \textit{split} or $(D',\Delta')$-\textit{split} if there is a number $D' \ge 0$ and a monotone function
$\Delta' \colon [0, + \infty ) \to [0, + \infty )$ so that
    \[
    F(U_1 \cup U_2) \subset F(U_1 [D',\Delta'_{x_0}]) + F(U_2 [D',\Delta'_{x_0}])
    \]
    for each pair of subsets $U_1$ and $U_2$ of $X \times Y$, \medskip
\item \textit{lean/split} or $(D,\Delta')$-\textit{lean/split} if there is a number $D \ge 0$ and a monotone function
$\Delta' \colon [0, + \infty ) \to [0, + \infty )$ so that \medskip
\begin{itemize}
\item the $X$-filtered module $F_X$ is $D$-lean, while
\item the $(X,Y)$-filtered module $F$ is $(D,\Delta')$-split,
\end{itemize} \medskip
\item \textit{insular} or $(d,\delta)$-\textit{insular} if there is a number $d \ge 0$ and a monotone function
$\delta \colon [0, + \infty ) \to [0, + \infty )$ so that
    \[
    F(U_1) \cap F(U_2) \subset F \big( U_1[d, \delta_{x_0}] \cap U_2[d, \delta_{x_0}] \big)
    \]
    for each pair of subsets $U_1$ and $U_2$ of $X \times Y$.
\end{itemize}
\end{DefRef}

\begin{PropRef}{HUNBV}
Suppose $F$ is an $(X,Y)$-filtered $R$-module.
\begin{enumerate}
\item If $F$ is $(D,\Delta)$-lean then the corresponding $X$-filtered module $F_X$ defined by assigning $F_X (S) = F(S \times Y)$ is $D$-lean.
\item Similarly, if $F$ is $(d,\delta)$-insular then $F_X$ is $d$-insular.
\item If $F$ is $(D,\Delta)$-lean then it is $(D,\Delta)$-split and, further, $(D,\Delta)$-lean/split.
\item An $(X,Y)$-filtered module $F$ which is lean/split and insular can be thought of as an object $F_X$ of $\LI (X,R)$.
\end{enumerate}
\end{PropRef}

\begin{proof}
	(1) Since $(x,y)[D,\Delta_{x_0}] \subset x[D] \times Y$, we have 
\[	F_X (S) \subset \sum_{x \in S} \sum_{y \in Y} F((x,y)[D,\Delta_{x_0}]) \subset \sum_{x \in S} F_X (x[D]). \]

(2) $F_X (S) \cap F_X (T) \subset F(S[d] \times Y \cap T[d] \times Y) = F_X (S[d] \cap T[d])$.

(3) The split property follows directly from definitions, and so the lean/split property follows in view of part (1).

(4) follows from (2).
\end{proof}

\begin{DefRef}{GEPCW3}
There are two subcategories nested in $\U_X (Y)$:
\begin{itemize}
\item $\mathbf{LS}_X (Y)$ is the full subcategory of $\U_X (Y)$ on objects $F$ that are lean/split and insular,
\item $\B_X (Y)$ is the full subcategory of $\mathbf{LS}_X (Y)$ on objects $F$ such that $F(U)$ is a finitely generated submodule whenever $U \subset X \times Y$ is bounded.  Equivalently, the subcategory $\B_X (Y)$ is full on objects $F$ such that all $Y$-filtered modules $F^S$ associated to bounded subsets $S \subset X$ are locally finitely generated.
\end{itemize}
\end{DefRef}

Clearly, $\B_X (Y)$ is a generalization of the bounded category $\B (X,R)$: if $Y = \point$ then $\B_X (Y)$ is precisely $\B (X,R)$.  On the other hand, if $X = \point$ then $\B_X (Y)$ is the full subcategory of $\BSI (Y,R)$ on locally finitely generated objects. 

We proceed to define appropriate exact structures in these categories.

\begin{DefRef}{ExUXY}
Let the
\textit{admissible monomorphisms} in $\U_X (Y)$ be the boundedly
bicontrolled homomorphisms $m \colon F_1 \to F_2$ such that the module homomorphism $F_1 (X \times Y) \to F_2 (X \times Y)$ is a monomorphism. Let the \textit{admissible epimorphisms} be the boundedly
bicontrolled homomorphisms $e \colon F_1 \to F_2$ such that $F_1 (X \times Y) \to F_2 (X \times Y)$ is an epimorphism.
The class $\mathcal{E}$ of
\textit{exact sequences} consists of the sequences
\[
F^{\subdot} \colon \ F' \xrightarrow{ \ i \ } F \xrightarrow{ \ j \ } F'',
\]
where $i$ is an admissible monomorphism, $j$ is an admissible epimorphism, and $\im (i) = \ker (j)$.
\end{DefRef}

\begin{ThmRef}{UbisWAb}
$\U_X (Y)$ is a Quillen exact category.
\end{ThmRef}

\begin{proof}
We will verify the axioms for exact structures due to Quillen with some simplifications due to B. Keller \cite{bK:90,bK:96}, cf. section 2 of \cite{bB:18}.

It follows from Lemma \refL{CatOK} that the collections of admissible monomorphisms and admissible epimorphisms are
closed under composition and that any short exact
sequence isomorphic to some sequence in
$\mathcal{E}$ is also in $\mathcal{E}$.

Now suppose we are given
an exact sequence $F' \xrightarrow{i } F \xrightarrow{j  } F''$
in $\mathcal{E}$ and a morphism $f \colon A \to F''$ in $\U_X
(Y)$.  Let $(b_j, \theta_j)$ be some filtration data for $j$ as a boundedly controlled epi and let $(b_f, \theta_f)$ be some contol data for $f$ as a boundedly controlled map.
There is a base change diagram
\[
\begin{CD}
F' @>>> F \times_{f} A @>{j'}>> A\\
@V{=}VV @VV{f'}V @VV{f}V \\
F' @>>> F @>{j}>> F''
\end{CD}
\]
where $m \colon F \times_{f} A \to F \oplus A$ is the kernel of
the epi
$j  \circ \pr_1 - f  \circ \pr_2 \colon F \oplus A \to F''$
and $f' =
\pr_1  \circ \, m$, $j' = \pr_2  \circ \, m$.
The $(X,Y)$-filtration on $F \times_{f} A$ is the standard filtration as a subobject of the product $F \times A$.
The induced map $j'$ has the same kernel as
$j$ and is bounded by $0$.
In fact,
\(
fA(U) \subset E''(U[b_f, \theta_{f,x_0}]),
\)
so
\(
fA(S) \subset j E(U[b_f + b_j, \theta_{f,x_0} + \theta_{j,x_0}]),
\)
and
\[
\im (j') \cap A(U) \subset j' \left( E \times_{f} A \right) \left(
U[b_f + b_j, \theta_{f,x_0} + \theta_{j,x_0}]] \right).
\]
This shows that $j'$ is boundedly
bicontrolled with filtration data $(b_f + b_j, \theta_{f} + \theta_{j})$.
Therefore, the class of admissible epimorphisms is
closed under base change by arbitrary morphisms in $\U_X
(Y)$.
Cobase changes by admissible monomorphisms are similar.
\end{proof}

\begin{PropRef}{SSSST}
The admissible monomorphisms are precisely the morphisms isomorphic in $\U_X (Y)$ to the filtration-wise monomorphisms and the admissible epimorphisms are those morphisms isomorphic to the filtration-wise epimorphisms.
In other words, the exact structure $\mathcal{E}$ in $\U_X
(Y)$ consists of sequences isomorphic to those
\[
E^{\subdot} \colon \quad E' \xrightarrow{\ i \ } E \xrightarrow{\ j \ } E''
\]
which possess filtration-wise restrictions
\[
E^{\subdot} (U) \colon \quad E' (U) \xrightarrow{\ i \ } E (U) \xrightarrow{\ j \ } E'' (U)
\]
for all subsets $U \subset (X,Y)$, and each $E^{\subdot} (U)$ is an exact sequence of $R$-modules.
\end{PropRef}

\begin{proof}
Each of the sequences $E^{\subdot}$ in the statement is an exact sequence in $\U_X (Y)$ because the restriction $i \colon E'(U) \to E(U)$ is monic and $j \colon E(U) \to E''(U)$ is epic, therefore $i \colon E' \to E$ and $j \colon E \to E''$ are both bicontrolled of filtration $(0,0)$.

Suppose $F^{\subdot}$ is a sequence isomorphic to such $E^{\subdot}$, so there is a commutative diagram
\[
\begin{CD}
F' @>{f}>> F @>{g}>> F''\\
@V{\cong}VV @VV{\cong}V @VV{\cong}V \\
E' @>{i}>> E @>{j}>> E''
\end{CD}
\]
Then $f$ and $g$ are compositions of two isomorphisms (which are clearly boundedly bicontrolled) which are either preceded by a boundedly bicontrolled monic or followed by a boundedly bicontrolled epi.  By part (1) of Lemma \refL{CatOK}, both  $f$ and $g$ are boundedly bicontrolled.

Now suppose $F^{\subdot}$ is an exact sequence in $\mathcal{E}$.  Let $K = \ker (g)$ and $C = \coker (f)$, then we obtain a commutative diagram
\[
\begin{CD}
F' @>{f}>> F @>{g}>> F''\\
@V{\cong}VV @VV{=}V @AA{\cong}A \\
K @>{i}>> F @>{j}>> C
\end{CD}
\]
where the vertical maps are the canonical isomorphisms.
By the construction of kernels and cokernels in the proof of Theorem \refT{MNBVX}, there are exact sequences
\(
K(U) \xrightarrow{ i  } F(U) \xrightarrow{ j  } C(U)
\)
for all subsets $U \subset (X,Y)$.
\end{proof}

\begin{PropRef}{TRGHJ1}
In the exact category $\U_X (Y)$,
\begin{enumerate}
\item the lean/split objects are closed under extensions,
\item the insular objects are closed under extensions.

\medskip
\noindent
Suppose
$E' \xrightarrow{f} E \xrightarrow{g} E''$
is an exact sequence in $\U_X (Y)$.

\medskip

\item If the object $E$ is
lean/split then $E''$ is lean/split.
\item If $E$ is
insular then $E'$ is insular.
\item Suppose $E$ is
insular then $E''$ is insular if $E'$ is lean/split.
\end{enumerate}
\end{PropRef}

\begin{proof}
All parts are proved by adapting the proofs of Lemmas \refL{LandS} and \refL{lninpres}.
To illustrate, 
suppose that in the exact sequence, as given in the statement, $(b, \theta)$ is common
filtration data for $f$ and $g$ and both $E'$ and $E''$ are
$(D,\Delta')$-lean/split.
For the first statement of part (2), notice that $E_X$ is $(4b+D)$-lean by part (1) of Lemma \refL{LandS}, so we need to verify that split objects are closed under extensions.
Consider two subsets $U_1$ and $U_2$ of $X \times Y$.  Then
\begin{align}
gE(U) &\subset E'' ((U_1 \cup U_2)[b, \theta_{x_0}]) \notag \\
&= E'' (U_1 [b, \theta_{x_0}] \cup U_2 [b, \theta_{x_0}]) \notag \\
&\subset E'' (U_1 [b + D, \theta_{x_0} + \Delta'_{x_0}]) + E'' (U_2 [b + D, \theta_{x_0} + \Delta'_{x_0}]). \notag
\end{align}
Therefore
\begin{align}
E(U) &\subset E (U_1 [2b + D, 2\theta_{x_0} + \Delta'_{x_0}]) + E (U_2 [2b + D, 2\theta_{x_0} + \Delta'_{x_0}])  \notag \\
&+ fE' (U_1 [3b + 2D, 3\theta_{x_0} + 2\Delta'_{x_0}]) + fE' (U_2 [3b + 2D, 3\theta_{x_0} + 2\Delta'_{x_0}]) \notag \\
&\subset E (U_1 [4b + 2D, 4\theta_{x_0} + 2\Delta'_{x_0}]) + E (U_2 [4b + 2D, 4\theta_{x_0} + 2\Delta'_{x_0}]),  \notag
\end{align}
showing that $E$ is $(4b+2D, 4\theta + 2\Delta')$-lean/split.
\end{proof}

\begin{ThmRef}{TRGHJ2}
$\mathbf{LS}_X (Y)$ is closed under extensions in $\U_X (Y)$.  In turn, $\B_X (Y)$  is closed under extensions in $\mathbf{LS}_X (Y)$.  Therefore, $\B_X (Y)$ is an exact category,
and the inclusion
$e \colon \mathcal{C}_X (Y) \to \B_X (Y)$
is an exact embedding.
\end{ThmRef}

\begin{proof}
The first statement follows from parts (2) and (3) of Proposition \refP{TRGHJ1}.
Suppose $f \colon F \to G$ is an isomorphism with $\fil (f) \le (b,\theta)$ and $G$ is locally finitely generated, then $F (U)$ is a finitely generated submodule of $G (U[b,\theta])$ for any bounded subset $U \subset X \times Y$ since $R$ is Noetherian.

If
\(
F' \xrightarrow{ f  } F \xrightarrow{ g  } F''
\)
is an exact sequence in $\mathbf{LS}_X (Y)$, $F'$ and $F''$ are locally finitely generated, and $(b, \theta)$ is common
filtration data for $f$ and $g$, then $gF(U)$ is a finitely generated submodule of $F'' (U[b,\theta])$ for any bounded subset $U$.
The kernel of the restriction of $g$ to $F(U)$ is a finitely generated submodule of $F'(U[b,\theta])$, so the extension $F(U)$ is finitely generated.
\end{proof}

\begin{RemRef}{TrKG}
(1) There is certainly an exact embedding
$\iota \colon \B (X \times Y, R) \to \B_X (Y)$
which is given by the identity on objects.
Because of the relaxation of the control conditions on homomorphisms, the morphism sets in the image of $\iota$ are in general properly smaller than in $\B_X (Y)$.
However $\iota$ is also proper on
objects. For example, the lean objects in $\BLI (X \times Y, R)$ are generated by the submodules $f (S \times T)$
where the diameters of $S$ and $T$ are uniformly bounded from above.
This is different from the weaker condition in
$\B_X (Y)$.

(2) While there is no functor between the categories $\U_X (Y)$ and $\U (X, \U (Y))$, there is a ``forgetful'' function associating to objects of $\U_X (Y)$ some objects of $\U (X, \U (Y))$.  It is defined by $\phi (F) (S) = F(S,Y)$ with the $Y$-filtration given by $F(S,Y) (T)= F(S,T)$.

This relationship can be made much more fruitful if the nature of the objects in $\U (X, \U (Y))$ is shifted to be functors from the category of \textit{bounded} subsets of $X$ to subobjects of $\U (Y)$ leading to a category that can be thought of as $\B (X, \B(Y,R))$.  It turns out that on this level there is a well-defined exact functor $\B_X (Y) \to \B (X, \B(Y,R))$.  We don't require this type of theory but it has been developed and applied in \cite{bB:18,pC:17}. 
\end{RemRef}

\SSecRef{Fibrewise restriction}{FibLoca}

We begin to prepare for the development of fibred localization exact sequences and fibred bounded excision theorems.
The model for localization and fibration theorems in controlled $G$-theory \cite[sections 3 and 4]{gCbG:00} can be implemented here as well. 

There are two complementary ways to introduce support in $\B_X (Y)$.
\begin{enumerate}
\item Let $\B_{<Z} (Y)$ be the full subcategory of $\B_X (Y)$ on objects $F$ supported near $Z$ when viewed as objects $F_X$ in $\BLI (X,R)$.
    In other words, $F$ is an object of $\B_{<Z} (Y)$ if
    \(
    F_X \subset F_X (Z[d]) = F(Z[d] \times Y)
    \)
    for some number $d \ge 0$.
\smallskip
\item Let $\B_X (Y)_{<C}$ be the full subcategory of $\B_X (Y)$ on objects $F$ such that
    \[
    F (X, Y) \subset F \big( (X, C)[r,\rho_{x_0}]) \big)
    \]
    for some number $r \ge 0$ and an order preserving function $\rho \colon [0,+\infty) \to [0,+\infty)$.
\end{enumerate}

The first version of support is a straightforward generalization of support for geometric modules that was exploited in \cite{gCbG:00}.
In this paper we are more interested in the latter, fibrewise version (2) of support.

\begin{PropRef}{HYUMN}
Suppose $F$ is a $(D,\Delta)$-lean/split object of $\B_X (Y)$.
The following are equivalent statements.
\begin{enumerate}
\item $F$ is an object of $\B_X (Y)_{<C}$.

\item There is a number $k \ge 0$ and an order preserving function $\lambda \colon \mathcal{B}(X) \to [0,+\infty)$ such that
    \[
    F^S \subset F^{S[k]}(C[\lambda(S)])
    \]
    for all bounded subsets $S \subset X$.
\item There is a number $k \ge 0$ and a monotone function $\Lambda \colon [0,+\infty) \to [0,+\infty)$
such that
\[
F^{x[D]} \subset F^{x [D+k]} (C[\Lambda_{x_0} (x)])
\]
for all $x \in X$.
\end{enumerate}
\end{PropRef}

\begin{proof}
$(2) \Longleftrightarrow (3)$: If $F$ satisfies (2) then
$F^{x[D]} \subset F^{x[D+k]} \big( C [\lambda (x[D])] \big)$.
It suffices to define $\Lambda$ such that
$\lambda(x[D]) \le \Lambda_{x_0} (x) = \Lambda (d(x_0,x))$.
Since $x[D] \le x_0 [d(x_0,x) +D]$ and $\lambda$ is order preserving,
one can take
$\Lambda (r) = \lambda \big( x_0 [r +D] \big)$.

In the opposite direction, given a bounded subset $S \subset X$,
\begin{equation}
F^S \subset\ \sum_{x \in S} F^{x[D]}
\subset\ \sum_{x \in S} F^{x [D+k]} (C[\Lambda_{x_0} (x)])
\subset\ F^{x [D+k]} (C[\lambda(S)])
\notag \end{equation}
when $\lambda (S) = \sup \{ \Lambda_{x_0} (x) \, \vert \, x \in S \}$.

$(1) \Longleftrightarrow (3)$: If $F$ is in $\B_X (Y)_{<C}$ then
$F^{x[D]} \subset F \big( (X, C) [r,\rho] \big)$,
so
\[
F^{x[D]} \subset F^{x[D]} \cap F \big( (X, C) [r,\rho] \big).
\]
If $F$ is $(d,\delta)$-insular then
\begin{equation}
F^{x[D]} \subset\ F \big( (x, C) [D+r+d, \rho + \delta] \big)
\subset\ F^{x [D+d +r]} (C[\Lambda_{x_0}(x)])
\notag \end{equation}
for $\Lambda (a) = \sup \{ (\delta + \rho)(z) \, \vert \, d(x_0,y) \le a +D+d+r \}$.

In the opposite direction, we have
\begin{equation}
F \subset\ \sum_{x \in X} F^{x[D]}
\subset\ \sum_{x \in X} F^{x [D+k]} (C[\Lambda_{x_0} (x)])
\subset\ F \big( (X, C) [D+k, \Lambda_{x_0}] \big)
\notag \end{equation}
for an object $F$ of $\B_X(Y)$ satisfying (3).
\end{proof}

\begin{DefRef}{GrCDef}
A \textit{Serre subcategory} of
an exact category is a full subcategory which is closed under
exact extensions and closed under passage to admissible subobjects and
admissible quotients.
\end{DefRef}

Note that this property is relative to the choice of exact structure. 

\begin{PropRef}{Serre3}
$\B_X (Y)_{<C}$ is a Serre subcategory of $\B_X (Y)$.
\end{PropRef}

\begin{proof}
First we show closure under exact extensions. Let
\(
F \xrightarrow{ f } G \xrightarrow{ g } H
\)
be an exact sequence in $\B_X (Y)$, let $(b, \theta)$ be common set of filtration
data for $f$ and $g$, and assume all objects be $(D,\Delta)$-lean/split.
We assume that $F$ and $H$ are objects of $\B_X (Y)_{<C}$, so there is a number $r \ge 0$ and a monotone function $\rho \colon [0,+\infty) \to [0,+\infty)$ such that
at the same time
$F (X, Y) = F \big( (X, C)[r,\rho_{x_0}]) \big)$ and
$H (X, Y) = H \big( (X, C)[r,\rho_{x_0}]) \big)$
for some choice of a base point $x_0$ in $X$.
Therefore
\begin{equation*}
fF(X, Y) = fF \big( (X, C)[r,\rho_{x_0}]) \big) \subset G \big( (X, C)[r+b,\rho_{x_0} +\theta_{x_0}]) \big).
\end{equation*}
In particular, the image $I = \im (f)$ with the standard filtration $I^S(T) = I \cap G^S(T)$ is an object of $\B_X (Y)_{<C}$.
Now
\begin{equation*}
H(X, Y) = gG (X,Y) \cap H \big( (X, C)[r,\rho_{x_0}]) \big) \subset gG \big( (X, C)[r+b,\rho_{x_0} +\theta_{x_0}]) \big).
\end{equation*}
Let $L = G \big( (X, C)[r+b,\rho_{x_0} +\theta_{x_0}] \big)$ viewed as a subobject of $G$ with the standard filtration.
Since $G = I + L$ for any submodule $L$ with $g(L)=H$, we have
\[
G (X, Y) = G \big( (X, C)[r+b,\rho_{x_0} +\theta_{x_0}]) \big),
\]
so $G$ is an object of $\B_X (Y)_{<C}$.

Suppose $f \colon F \to G$ is an admissible monomorphism in $\B_X (Y)$, which is a boundedly bicontrolled monic with $\fil (f) \le
(b,\theta)$, $F$ is $(D',\Delta')$-lean/split, $G$ is $(D,\Delta)$-lean/split for some $D \ge D' + b$, and
$G$ is $(d,\delta)$-insular.

If $G$ is an object of $\B_X (Y)_{<C}$, according to Proposition \refP{HYUMN},
\[
G^S \subset G^{S[k]}(C[\lambda(S)])
\]
for some number $k \ge 0$, an order preserving function $\lambda \colon \mathcal{B} (X) \to [0, +\infty)$, and all bounded subsets $S \subset X$.
Then
\begin{multline}
f F^{x[D']} \subset G^{x[D' + b]} \subset G^{x[D' + b + k]} (C[\lambda (x[D' + b])]) \subset G^{x[D + k]} (C[\lambda (x[D])]),
\notag
\end{multline}
using the fact that $\lambda$ is order preserving.
Since
\[
G^{x[D + k]} (Y - C[\lambda(x[D+k]) + \Delta(x[D]) + \theta(x[D']) + 2\delta(x[D+k])] ) = 0,
\]
we have
\[
F^{x[D']} (Y - C[\lambda(x[D+k]) + \Delta(x[D]) + 2\theta(x[D']) + 2\delta(x[D+k])] ) = 0.
\]
Therefore
\[
F^{x[D']} \subset F^{x[D']} (C[\lambda(x[D]) + \Delta(x[D]) + \Delta'(x[D]) + 2\theta(x[D']) + 2\delta(x[D+k])],
\]
so $F$, which is generated by $F^{x[D']}$, is also an object of $\B_X (Y)_{<C}$.

On the
other hand, let $g \colon G \to H$ be an admissible quotient with
$\fil (g) \le (b, \theta)$
and suppose $G$ is an object of $\B_X (Y)_{<C}$ so that there is a number $r \ge 0$ and a monotone function $\rho \colon [0,+\infty) \to [0,+\infty)$ such that
\[
G (X, Y) = G \big( (X, C)[r,\rho_{x_0}]) \big).
\]
This implies that
\begin{equation*}
H(X, Y) = gG (X, Y) \subset H \big( (X, C)[r+b,\rho_{x_0} +\theta_{x_0}]) \big),
\end{equation*}
so
$H$ is also in $\B_X (Y)_{<C}$.
\end{proof}


\SSecRef{Fibrewise gradings}{FFFGGG}

The gradings from Definition \refD{WEI} can be generalized to gradings of objects from $\B_X (Y)$. 

\begin{DefRef}{LScov2}
Given an object $F$ of $\B_X (Y)$, a \textit{grading} of $F$
is a covariant functor
$\mathcal{F} \colon \mathcal{P}(X,Y) \to \mathcal{I}(F)$
with the following properties:
\begin{enumerate}
\item if $\mathcal{F} (\mathcal{C})$ is given the standard filtration, it is an object of $\B_X (Y)$,
\item there is an enlargement data $(K,k)$ such that
\[
F(C) \subset \mathcal{F} (C) \subset F(C [K,k_{x_0}]),
\]
for all subsets $C$ of $(X,Y)$.
\end{enumerate}
\end{DefRef}

\begin{RemRef}{OINBFGY2}
If $C = (X,S)$ then $\mathcal{F} (C)$ is an object of $\B_X (Y)_{<S}$.
\end{RemRef}

We are concerned with localizations to a specific type of subspaces of $(X,Y)$.
This makes the following partial gradings sufficient and easier to work with.

\begin{DefRef}{LScov22}
Let $\mathcal{M}^{\ge 0}$ be the set of all monotone functions $\delta \colon [0, + \infty) \to [0, + \infty)$.
Let $\mathcal{P}_X (Y)$ be the subcategory of $\mathcal{P}(X,Y)$ consisting of
all subsets of the form $(X,C)[D,\delta_{x_0}]$ for some choices of a subset $C \subset Y$, a number $D \ge 0$, and a function $\delta \in \mathcal{M}^{\ge 0}$.

Given an object $F$ of $\B_X (Y)$, a $Y$-\textit{grading} of $F$
is a functor
$\mathcal{F} \colon \mathcal{P}_X (Y) \to  \mathcal{I}(F)$
with the following properties:
\begin{enumerate}
\item the submodule $\mathcal{F} ((X,C)[D,\delta_{x_0}])$ with the standard filtration is an object of $\B_X (Y)$,
\item there is an enlargement data $(K,k)$ such that
\[
F((X,C)[D,\delta_{x_0}]) \subset \mathcal{F} ((X,C)[D,\delta_{x_0}]) \subset F((X,C) [D + K, \delta_{x_0} + k_{x_0}]),
\]
for all subsets in $\mathcal{P}_X (Y)$.
\end{enumerate}

Since $U [D + K, \delta_{x_0} + k_{x_0}] = U [D, \delta_{x_0}] [K, k_{x_0}]$ for general subsets $U$, the third, largest submodule is independent of the choice of $D$, $\delta_{x_0}$.
\end{DefRef}

We say that an object $F$ of $\B_X (Y)$ is $Y$-\textit{graded} if there exists a $Y$-grading of $F$, but the grading itself is not specified, and define $\G_X (Y)$ as the full subcategory of $\B_X (Y)$ on $Y$-graded filtered modules.

\begin{PropRef}{VFGHJ2}
The $Y$-graded objects in $\B_X (Y)$ are closed under isomorphisms.
The subcategory $\G_X (Y)$ is closed under extensions in $\B_X (Y)$.
Therefore, $\G_X (Y)$ is an exact subcategory of $\B_X (Y)$.
\end{PropRef}

\begin{proof}
Suppose $f \colon F \to F'$ is an isomorphism between a $Y$-graded module $F$ and $F'$ in $\B_X (Y)$.  If $\fil (f) \le (b,\theta)$ and $\mathcal{F}$ is a $Y$-grading for $F$ then $\mathcal{F}' (X,C)[D,\delta_{x_0}] = f F((X,C) [D + K + b, \delta_{x_0} + k_{x_0} + \theta_{x_0}])$.  The rest of the argument closely follows the proof of Proposition \refP{VFGHJ}.  We want to spell out one detail for future reference.  Let 
\[
F \xrightarrow{\ f \ } G \xrightarrow{\ g \ } H
\]
be an exact sequence in $\B_X (Y)$.  Suppose $\fil (f) \le (b,\theta)$ and $\fil (g) \le (b,\theta)$ for the same set of bicontrol bound data and suppose that
$F$ and $H$ are graded modules in $\G_X (Y)$ with the associated functors $\mathcal{F}$ and $\mathcal{H}$.
The assignment
\[
\mathcal{G}(U) = f \mathcal{F} (U[3b,3\theta]) + G(S[2b,2\theta]) \cap g^{-1} \mathcal{H} (S[b,\theta])
\]
gives a grading of $G$ as an object of $\G_X (Y)$.
\end{proof}

As with the category $\G (X,R)$, the advantage of working with $\G_X (Y)$ as opposed to $\B_X (Y)$ is that we are able to localize to the grading subobjects associated to subsets from the family $\mathcal{P}_X (Y)$.

\begin{LemRef}{JHQASF22}
Let $F$ be a submodule of a $Y$-filtered module $G$ in $\G_X (Y)$ which is lean/split with respect to the standard filtration.  Then $\mathcal{F} (U) = F \cap \mathcal{G} (U)$ is a $Y$-grading of $F$.
\end{LemRef}

We will call this induced $Y$-grading of $F$ the \textit{standard $Y$-grading} of the submodule.

\begin{proof}
The proof is reduced to checking that $\mathcal{F}(U)$ is an object of $\B_X (Y)$ for each subset $U \in \mathcal{P}_X (Y)$.
Suppose $i \colon F \to G$ is the inclusion and $q \colon G \to H$ is the quotient of $i$.
Since $F$ is insular by part (4) of Proposition \refP{TRGHJ1}, both $F$ and $G$ are lean/split and insular.
Thus $H$ is lean/split and insular by parts (3) and (5) of \refP{TRGHJ1}.
Let $\mathcal{H}(U) = q \mathcal{G}(U)$ with the standard filtration in $H$. Then $\mathcal{H}(U)$ is lean/split by part (3) and insular as a submodule of insular $H$.
The kernel $\mathcal{F}(U)$ of the filtration $(0,0)$ map $q \vert \colon \mathcal{G}(U) \to \mathcal{H}(U)$
is lean/split by part (5) of \refP{TRGHJ1} and is insular as a submodule of insular $F$.
Locally finite generation of $\mathcal{F}(U)$ follows from the same property of $\mathcal{G}(U)$.
\end{proof}

\begin{PropRef}{HUVZFMNtt}
Suppose $g \colon G \to H$ is a boundedly bicontrolled epimorphism in $\B_X (Y)$ and suppose $F$ is the kernel of $g$ in $\Mod (R)$.  If $G$ is $Y$-graded and $F$ is lean/split with respect to the standard $Y$-filtration then both $H$ and $F$ are $Y$-graded.
\end{PropRef}

\begin{proof}
The $Y$-grading for $H$ is given by $\mathcal{H}(U) = g \mathcal{G}(U[b,\theta])$, where $(b, \theta)$ is a chosen set of filtration data for $g$.
The argument for Lemma \refL{JHQASF22} applies directly to show that $\mathcal{H}$ is indeed a grading of $H$, and the assignment $\mathcal{F}(U) = F \cap \mathcal{G}(U[b,\theta])$ gives a $Y$-grading of $F$.
\end{proof}

This allows to characterize admissible monomorphisms in $\G_X (Y)$.

\begin{PropRef}{NJWSXnw}
The inclusion of a subobject $i \colon F \to G$ in $\G_X (Y)$ is an admissible monomorphism if and only if $F$ is lean/split.	
\end{PropRef}

\begin{proof}
Let $H$ be a cokernel of $i$ in $\U_X (Y)$.  As verified in Lemma \refL{JHQASF22}, $H$ is lean/split and insular and is, in fact, a cokernel of $i$ in $\B_X (Y)$. 
From Proposition \refP{HUVZFMNtt}, $H$ is graded, so it is also a cokernel of $i$ in $\G(X,R)$.
\end{proof}

In the case $K \le 0$ and $k$ is a non-positive function, we will use notation $U [K,k]$ for the subset $U \setminus ((X,Y) \setminus U) [-K,-k]$ of $(X,Y)$.  The following is a direct analogue of Corollary \refC{LPCSL} with exactly the same proof. 

\begin{CorRef}{LPCSLtt}
Given an object $F$ in $\G_X (Y)$ and a subset $U$ from the family $\mathcal{P}_X (Y)$, there is a set of enlargement data $(K,k)$ and an admissible subobject $i \colon F_U \to F$ in $\G_X (Y)$ with the property that $F_U \subset F(U[K,k])$.
If $G$ is $(D,\Delta')$-lean/split then the quotient $q \colon F \to H$ of the inclusion has the property that $H (X,Y) = H (((X,Y) \setminus U) [2D,2\Delta'])$.
\end{CorRef}

Now we can summarize the preceding results.

\begin{ThmRef}{LPCSLuu}
Given a graded object $F$ in $\G_X (Y)$ and a subset $U$ from the family $\mathcal{P}_X (Y)$,
we assume that $F$ is $(D,\Delta')$-split and $(d, \delta)$-insular and is graded by $\mathcal{F}$.
The submodule $\mathcal{F} (U)$ has the following properties:
\begin{enumerate}
\item $\mathcal{F} (U)$ is graded by $\mathcal{F}_U (T) = \mathcal{F} (U) \cap \mathcal{F} (T)$,
\item $F(U) \subset \mathcal{F} (U) \subset F(U[K,k])$ for some fixed enlargement data $(K,k)$,
\item if $q \colon F \to H$ is the quotient of the inclusion $i \colon \mathcal{F} (U) \to F$ and $F$ is $(D, \Delta')$-lean/split, then $H$ is supported on $(X \setminus U) [2D,2\Delta']$,
\item $H (U[-2D -2d, -2\Delta' - 2\delta]) =0$.
\end{enumerate}
\end{ThmRef}

\begin{proof}
Property (1) follows from Lemma \refL{JHQASF22}. 
Properties (2) and (3) follow from Corollary \refC{LPCSLtt}. 
(4) follows from the fact that a $d$-insular filtered module is $(2d,2\delta)$-separated,
in the sense that for any pair of subsets $U$ and $V$ of $(X,Y)$ such that $U[2d, 2\delta] \cap V = \emptyset$ we have
$U[d,\delta] \cap V[d,\delta] = \emptyset$ so $F(U) \cap F(V) = 0$.
Now 
\[
H (U[-2D -2d, -2\Delta' - 2\delta]) \cap H \left( ((X,Y) \setminus U)[2D,2\Delta'] \right) =0, 
\]
but
$H(((X,Y) \setminus U)[2D,2\Delta']) = H(X,Y)$, thus $H (U[-2D -2d, -2\Delta' - 2\delta]) = 0$.
\end{proof}

\SSecRef{Localization fibration sequence}{LLFFSS}

We will use the localization theorem of Schlichting \cite{mS:03} for Serre subcategories of exact categories.  These techniques require the Serre subcategory to satisfy some additional assumptions that we verify next.

\begin{DefRef}{CLFracs}
A class of morphisms $\Sigma$ in an additive category $\mathcal{A}$
\textit{admits a calculus of right fractions} if
\begin{enumerate}
\item the identity of each object is in $\Sigma$,
\item $\Sigma$ is closed under composition,
\item each diagram $F \xrightarrow{\ f} G \xleftarrow{\ s\ } G'$ with $s
\in \Sigma$ can be completed to a commutative square
\[
\begin{CD}
 F' @>{f'}>> G'\\
 @VV{t}V @VV{s}V\\
 F @>{f}>> G
\end{CD}
\]
with $t \in \Sigma$, and
\item if $f$ is a morphism in $\mathcal{A}$ and $s \in \Sigma$ such that
$sf = 0$ then there exists $t \in \Sigma$ such that $ft = 0$.
\end{enumerate}
In this case there is a construction of the \textit{localization}
$\mathcal{A} [\Sigma^{-1}]$ which has the same objects as $\mathcal{A}$.  The
morphism sets $\Hom (F,G)$ in $\mathcal{A} [\Sigma^{-1}]$ consist of
equivalence classes of diagrams
\[
(s,f) \colon \quad F \xleftarrow{\ s\ } F' \xrightarrow{\ f} G
\]
with the equivalence relation generated by $(s_1,f_1) \sim
(s_2,f_2)$ if there is a map $h \colon F'_1 \to F'_2$ so that $f_1
= f_2 h$ and $s_1 = s_2 h$. Let $(s \vert f)$ denote the
equivalence class of $(s,f)$. The composition of morphisms in $\mathcal{A}
[\Sigma^{-1}]$ is defined by
$(s \vert f) \circ (t \vert g) = (st' \vert gf')$
where $f'$ and $t'$ fit in the commutative square
\[
\begin{CD}
 F'' @>{f'}>> G'\\
 @VV{t'}V @VV{t}V\\
 F @>{f}>> G
\end{CD}
\]
from axiom 3.
\end{DefRef}

\begin{PropRef}{FactsFrac}
The localization $\mathcal{A} [\Sigma^{-1}]$
is a category.  The morphisms
of the form $(\id \vert s)$ where
$s \in \Sigma$ are isomorphisms
in $\mathcal{A} [\Sigma^{-1}]$. The rule
$P_{\Sigma} (f) = (\id \vert f)$
gives a functor
$P_{\Sigma} \colon \mathcal{A} \to
\mathcal{A} [\Sigma^{-1}]$ which
is universal among the functors making
the morphisms $\Sigma$
invertible.
\end{PropRef}

\begin{proof}
The proofs of these facts can be found in Chapter I of
\cite{pGmZ:67}.  The inverse of $(\id \vert s)$ is $(s \vert
\id)$.
\end{proof}

From Proposition \refP{Serre3}, given a subset $C$ of $Y$, the category $\B_X (Y)_{<C}$ is a Serre subcategory of $\B_X (Y)$.

\begin{PropRef}{NMCWJK}
The restriction to $Y$-gradings in $\B_X (Y)_{<C}$ gives a full exact subcategory $\G_X (Y)_{<C}$ which is a Serre subcategory of $\G_X (Y)$.
\end{PropRef}

\begin{proof}
One only needs to observe that  if the modules at the ends of the exact sequence in the proof of Proposition \refP{VFGHJ2} have $Y$-gradings in $\B_X (Y)_{<C}$  then the displayed grading of $G$ shows that $G$ is an object of $\G_X (Y)_{<C}$.
\end{proof}

The following shorthand notation is convenient when the choice of $C$ is clear.

\begin{NotRef}{HUY}
The category $\G$ is the exact subcategory of $Y$-graded objects in $\B_X (Y)$.
When the choice of the subset $C \subset Y$ is understood, we will use notation $\C$ for the Serre subcategory
$\G_X (Y)_{<C}$ of $\G$.
\end{NotRef}

\begin{DefRef}{YTRE}
Define the class
of \textit{weak equivalences} $\Sigma (C)$ in $\G$ to consist of all finite
compositions of admissible monomorphisms with cokernels in $\C$ and
admissible epimorphisms with kernels in $\C$.
\end{DefRef}

We need the class $\Sigma (C)$ to admits calculus of right fractions.
This follows from \cite[Lemma 1.13]{mS:03} as soon as we prove the following fact.

A Serre subcategory $\C$ of an exact category $\G$ is \textit{right
filtering} if each morphism $f \colon F_1 \to F_2$ in $\G$, where $F_2$
is an object of $\C$, factors through an admissible epimorphism $e
\colon F_1 \to \overline{F}_2$, where $\overline{F}_2$ is in $\C$.

\begin{LemRef}{SubSubLem1}
The subcategory $\C = \G_X (Y)_{<C}$ of $\G = \G_X (Y)$ is right filtering.
\end{LemRef}

\begin{proof}
For a morphism $f \colon F_1 \to F_2$ in $\G$ with $F_2$
in $\C$, we assume that both $F_1$ and $F_2$ are
$(D, \Delta')$-lean/split and $(d, \delta)$-insular.
Suppose $f$ is bounded by $(b,\theta)$ and let $r \ge 0$ and $\rho \colon [0,+\infty) \to [0,+\infty)$ be a
monotone function such that
\[
F_2 (X,Y) \subset F_2 ((X,C)[r, \rho_{x_0}]).
\]
Now for any characteristic set of data $(K,k)$ for the grading $\mathcal{F}_1$ and any subset $R$ we have
\[
f \mathcal{F}(R) \subset fF_1 (R[K,k_{x_0}]) \subset F_2 (R[K + b, k_{x_0} + \theta_{x_0}]).
\]
By part (3) of Theorem \refT{LPCSLuu},
$F_2 (R[K + b, k_{x_0} + \theta_{x_0}]) \cap F_2 ((X,C)[r, \rho_{x_0}]) = 0$
for any $R$ such that
\[
R[K + b + 2D + 2d, k_{x_0} + \theta_{x_0} + 2\Delta'_{x_0} + 2\delta_{x_0}] \cap (X,C)[r, \rho_{x_0}] = \emptyset.
\]
If we choose
\[
R = (X,Y) \setminus (X,C) [K + b + 2D + 2d + r, k_{x_0} + \theta_{x_0} + 2\Delta'_{x_0} + 2\delta_{x_0} + \rho_{x_0}]
\]
and define $E = \mathcal{F}_1 (R)$, then
$f E = 0$.
Let $\overline{F}_2$ be the cokernel of the inclusion $E \to F_1$.
Then $\overline{F}_2$ is lean/split and insular and has a grading given by $\overline{\mathcal{F}}_2 (S) = q \mathcal{F}_1 (S[b,\theta_{x_0}])$.  Since
\[
\overline{F}_2 (X,Y) \subset \overline{F}_2 ((X,C) [K + b + 2D + 2d + r, k_{x_0} + \theta_{x_0} + 2\Delta'_{x_0} + 2\delta_{x_0} + \rho_{x_0}]),
\]
the quotient $\overline{F}_2$ is in $\C$, and $f$ factors as $F_1 \to \overline{F}_2 \to F_2$ in the right square in the map of exact sequences
\[
\begin{CD}
E @>>> F_1 @>{j'}>> \overline{F}_2\\
@V{i}VV @VV{=}V @VVV \\
K @>{k}>> F_1 @>{f}>> F_2
\end{CD}
\]
as required.
\end{proof}

\begin{DefRef}{Quot}
The category $\mathbf{G}/\mathbf{C}$ is the localization $\G [\Sigma (C)^{-1}]$.
\end{DefRef}

It is clear that the quotient $\mathbf{G}/\mathbf{C}$ is an
additive category, and $P_{\Sigma (C)}$ is an additive functor.
In fact, we have the following.

\begin{ThmRef}{ExLocStr}
The short sequences in $\mathbf{G}/\mathbf{C}$ which are
isomorphic to images of exact sequences from $\G$ form a Quillen
exact structure.
\end{ThmRef}

\begin{proof}
This will be a consequence from \cite[Proposition 1.16]{mS:03}.
Since $\C$ is right filtering by Lemma \refL{SubSubLem1}, it remains
to check that $\C$ right s-filtering in $\G$ in the following sense.
A subcategory $\C$ of an exact category $\G$ is \textit{right s-filtering} if given an admissible
monomorphism $f \colon F_1 \to F_2$ with $F_1$ in $\C$, there exist $E$ in $\C$ and an
admissible epimorphism $e \colon F_2 \to E$ such that the composition $ef$ is an
admissible monomorphism.

Suppose that $F_1$ and $F_2$ have the same properties as in the proof of Lemma \refL{SubSubLem1}, and
$\fil (f) \le (b, \theta)$.
Since $F_1$ is in $\C$, there are $r \ge 0$ and a monotone function
$\rho \colon [0,+\infty) \to [0,+\infty)$ such that
$F_1 (X,Y) \subset F_1 \left( (X,C) [r,\rho_{x_0}] \right)$.
Then let $F'_2 = \mathcal{F}_2 (T)$ where
\[
T = (X,Y) \setminus (X,C) [K + b + 2D + 2d + r, k_{x_0} + \theta_{x_0} + 2\Delta'_{x_0} + 2\delta_{x_0} + \rho_{x_0}].
\]
Define $E$ as the cokernel of the inclusion $F'_2 \to F_2$ and let $e \colon F_2 \to E$ be the quotient map.
The composition $ef$ is an admissible
monomorphism with $\fil (ef) = \fil (f) \le (b, \theta)$.
\end{proof}

\begin{NotRef}{GHPL}
If $C$ is a subset of $Y$ as before, $\G_X (Y,C)$ will stand for the exact category ${\G}/{\C}$ and $G_X (Y,C)$ for its Quillen $K$-theory.
\end{NotRef}

The main tool in proving controlled
excision theorems will be the following localization sequence.

\begin{ThmRefName}{Schlichting}{Theorem 2.1 of Schlichting \cite{mS:03}}
Let $\Z$ be an idempotent complete right s-filtering subcategory
of an exact category $\E$ which is full and closed under exact extensions.  Then the sequence of exact categories
$\Z \to \E \to {\E}/{\Z}$ induces a homotopy fibration of Quillen $K$-theory
spectra
\[
K(\Z) \longrightarrow K(\E) \longrightarrow K({\E}/{\Z}).
\]
\end{ThmRefName}

\begin{CorRef}{FibrG}
There is a homotopy fibration
\begin{gather*}
G_X (Y)_{<C} \longrightarrow G_X (Y) \longrightarrow G_X (Y,C).
\end{gather*}
\end{CorRef}

There is a more intrinsic statement of the same fact.

\begin{ThmRefName}{LocText}{Localization}
There is a homotopy fibration
\begin{gather*}
{G}_X (C) \longrightarrow {G}_X (Y) \longrightarrow {G}_X (Y,C).
\end{gather*}
\end{ThmRefName}

Theorem \refT{LocText} is a consequence of Corollary \refC{FibrG}
as soon as we show that ${G}_X (C)$ and ${G}_X (Y)_{<C}$ are weakly equivalent.

Recall that the \textit{essential full image} of a functor
$F \colon \C \to \D$ is the full subcategory of $\D$ whose objects
are those $D$ that are isomorphic to $F(C)$ for some $C$ from $\C$.

\begin{LemRef}{LemNeed1}
Given a pair of proper metric spaces $C \subset Y$, there is a fully faithful embedding $\epsilon \colon \G_X (C) \to
\G_X (Y)$.
The Serre subcategory $\G_X (Y)_{<C}$ is the essential full image of $\G_X (C)$ in $\G_X (Y)$.
Therefore, the inclusion $C \subset Y$ induces
a weak equivalence
\begin{gather*}
{G}_X (C) \longrightarrow {G}_X (Y)_{<C}.
\end{gather*}
\end{LemRef}

\begin{proof}
Suppose $F$ is an object of $\G_X (C)$.
The embedding $\epsilon$ is given by $\epsilon (F)(U) = F ((X,C) \cap U)$, $\epsilon (\mathcal{F})(S) = \mathcal{F} ((X,C) \cap U)$.
It is clear that $\epsilon (F)$ is in $\G_X (Y)_{<C}$.

To show that $\G_X (Y)_{<C}$ is the essential full image, for an object $G$ of $\G_X (Y)_{<C}$ assume that
$G \subset G((X,C)[r,\rho_{x_0}])$
for some number $r \ge 0$ and a monotone function $\rho \colon [0,+\infty) \to [0,+\infty)$.
Choose any set function
$\tau \colon (X,C)[r,\rho_{x_0}] \rightarrow X \times C$
with the properties
\begin{enumerate}
\item $\tau (x,y) = (x, \tau_x (y))$,
\item $d(y,\tau_x (y)) \le \rho_{x_0} +r$ for all $x$ in $X$,
\item $\tau \vert (X \times C) = \id$.
\end{enumerate}
Then the $Y$-filtered module $E$ associated to $G$ given by $E(S) = G (\tau^{-1}(S))$ with the grading
$\mathcal{E} (U') = \mathcal{G} (\tau^{-1} (U'))$ is an object of $\G_X (C)$.
Indeed, if $\mathcal{G} (\tau^{-1} (U'))$ is $(D,\Delta')$-lean/split and $(d, \delta)$-insular then $\mathcal{E} (U')$ is $(D + r,\Delta' + \rho)$-lean/split and $(d + r, \delta + \rho)$-insular.
The identity map is an isomorphism in $\G_X (Y)$ with $\fil (\id) \le (2r, 2\rho + 2r)$.
\end{proof}

\SecRef{Fibrewise excision theorems}{FWBET}

\SSecRef{Waldhausen categories and \textit{K}-theory}{WWW}

Our main reference for Waldhausen $K$-theory terminology and notation is Thomason \cite{rTtT:90}. 

A Waldhausen category $\D$ with weak equivalences $\w (\D)$ is
often denoted by $\wD$ as a reminder of the choice. A functor
between Waldhausen categories is exact if it preserves the chosen zero objects,
cofibrations, weak equivalences, and cobase changes.
Let $\D$ be a small Waldhausen category with respect to two
categories of weak equivalences $\bsv(\D) \subset \w(\D)$ with a
cylinder functor $T$ both for $\vD$ and for $\wD$ satisfying the
cylinder axiom for $\wD$. Suppose also that $\w(\D)$ satisfies the
extension and saturation axioms. Define $\vDw$ to be the full
subcategory of $\vD$ whose objects are $F$ such that $0 \to F \in
\w(\D)$. Then $\vDw$ is a small Waldhausen category with
cofibrations $\co (\Dw) = \co (\D) \cap \Dw$ and weak
equivalences $\bsv (\Dw) = \bsv (\D) \cap \Dw$. The cylinder
functor $T$ for $\vD$ induces a cylinder functor for $\vDw$.  If
$T$ satisfies the cylinder axiom then the induced functor does so
too.

\begin{ThmRefName}{ApprThm}{Approximation Theorem}
Let $E \colon \D_1 \rightarrow \D_2$ be an exact functor between
two small saturated Waldhausen categories. It induces a map of
$K$-theory spectra
\[
  K (E) \colon K (\D_1) \longrightarrow K (\D_2).
\]
Assume that $\D_1$ has a cylinder functor satisfying the cylinder
axiom. If $E$ satisfies two conditions:
\begin{enumerate}
\item a morphism $f \in \D_1$ is in $\w(\D_1)$
if and only if $E (f) \in \D_2$ is in $\w(\D_2)$,
\item for any object $D_1 \in \D_1$ and any morphism
$g \colon E(D_1) \to D_2$ in $\D_2$, there is an object $D'_1 \in
\D_1$, a morphism $f \colon D_1 \to D'_1$ in $\D_1$, and a weak
equivalence $g' \colon E(D'_1) \rightarrow D_2 \in \w(\D_2)$ such
that $g = g' E(f)$,
\end{enumerate}
then $K (E)$ is a homotopy equivalence.
\end{ThmRefName}

\begin{proof}
This is Theorem~1.6.7 of \cite{fW:85}. The presence of the
cylinder functor with the cylinder axiom allows to make condition (2)
weaker than that of Waldhausen, see point 1.9.1 in
\cite{rTtT:90}.
\end{proof}

\begin{DefRef}{PFff}
In any additive category, a sequence of morphisms
\[
E^{\subdot} \colon \quad 0 \longrightarrow E^1 \xrightarrow{\ d_1
\ } E^2 \xrightarrow{\ d_2 \ }\ \dots\ \xrightarrow{\ d_{n-1} \ }
E^n \longrightarrow 0
\]
is called a \textit{(bounded) chain complex} if the compositions
$d_{i+1} d_i$ are the zero maps for all $i = 1$,\dots, $n-1$.  A
\textit{chain map} $f \colon F^{\subdot} \to E^{\subdot}$ is a
collection of morphisms $f^i \colon F^i \to E^i$ such that $f^i
d_i = d_i f^i$.  A chain map $f$ is \textit{null-homotopic} if
there are morphisms $s_i \colon F^{i+1} \to E^i$ such that $f = ds
+ sd$.  Two chain maps $f$, $g \colon F^{\subdot} \to E^{\subdot}$
are \textit{chain homotopic} if $f-g$ is null-homotopic. Now $f$
is a \textit{chain homotopy equivalence} if there is a chain map
$h \colon E^i \to F^i$ such that the compositions $fh$ and $hf$
are chain homotopic to the respective
identity maps.
\end{DefRef}

The Waldhausen structures on categories of bounded chain complexes
are based on homotopy equivalence as a weakening of the notion of
isomorphism of chain complexes.

A sequence of maps in an exact category is called \textit{acyclic}
if it is assembled out of short exact sequences in the sense that
each map factors as the composition of the cokernel of the
preceding map and the kernel of the succeeding map.

It is known that the class of acyclic complexes in an exact
category is closed under isomorphisms in the homotopy category if
and only if the category is idempotent complete, which is also
equivalent to the property that each contractible chain complex is
acyclic, cf.\ \cite[sec.\ 11]{bK:96}.

Given an exact category $\E$, there is a standard choice for the
Waldhausen structure on the category $\E'$ of bounded
chain complexes in $\E$ where the degree-wise admissible
monomorphisms are the cofibrations and the chain maps whose
mapping cones are homotopy equivalent to acyclic complexes are the
weak equivalences $\boldsymbol{v}(\E')$.

The following fact is well-known, cf. point
1.1.2 in \cite{rTtT:90}.

\begin{PropRef}{FibThApp2}
The category $\vE'$ is a Waldhausen category satisfying the
extension and saturation axioms and has cylinder functor
satisfying the cylinder axiom.
\end{PropRef}

\begin{ExRef}{FibSt77}
There are two choices for the Waldhausen structure on the category
of bounded chain complexes $\G' = \G (X,R)'$. One is the standard choice $\vG'$ as above. Given a subset $C \subset Y$,
another choice for the weak equivalences $\w (\G')$ is the chain
maps whose mapping cones are homotopy equivalent to acyclic
complexes in the quotient $\mathbf{G}/\mathbf{C}$.
\end{ExRef}

\begin{CorRef}{FibSt}
	The categories $\vG'$ and $\wG'$ are Waldhausen categories
satisfying the extension and saturation axioms and have cylinder
functors satisfying the cylinder axiom.
\end{CorRef}

\begin{proof}
	All axioms and constructions, including the cylinder functor, for
$\wG'$ are inherited from $\vG'$.
\end{proof}

The $K$-theory functor from the category of small Waldhausen
categories $\D$ and exact functors to the category of connective spectra is
defined in terms of $S_{\subdot}$-construction as in Waldhausen
\cite{fW:85}. It extends to simplicial categories $\D$ with
cofibrations and weak equivalences and inductively delivers the
connective spectrum $n \mapsto \vert \bfw S_{\subdot}^{(n)} \D
\vert$. We obtain the functor assigning to $\D$ the connective
$\Omega$-spectrum
\[
K (\D) = \Omega^{\infty} \vert \bfw S_{\subdot}^{(\infty)} \D
\vert = \colim{n \ge 1} \Omega^n \vert \bfw S_{\subdot}^{(n)} \D
\vert
\]
representing the Waldhausen algebraic $K$-theory of $\D$. For
example, if $\D$ is the additive category of free finitely
generated $R$-modules with the canonical Waldhausen structure,
then the stable homotopy groups of $K (\D)$ are the usual
$K$-groups of the ring $R$.  In fact, there is a general
identification of the two theories. Recall that for any exact
category $\E$, the category $\E'$ of bounded chain complexes has the Waldhausen
structure $\vE'$ as in Example \refE{FibSt77}.

\begin{ThmRef}{Same}
The Quillen $K$-theory of an exact category $\E$ is equivalent to
the Waldhausen $K$-theory of $\vE'$.
\end{ThmRef}

\begin{proof}
The proof is based on repeated applications of the Additivity
Theorem, cf. Thomason's Theorem 1.11.7 from \cite{rTtT:90}. Thomason's
proof of his Theorem 1.11.7 can be repeated verbatim here.  It is
in fact simpler in this case since his condition 1.11.3.1 is not
required.
\end{proof}

\SSecRef{Controlled excision theorems}{ET}

These are the major computational tools in controlled $K$-theory.
We develop excision results $G_X (Y)$ with respect to specific coverings of the variable $Y$ in this section.  

Suppose $Y_1$
and $Y_2$ are subsets of a proper metric space $Y$, and $Y = Y_1
\cup Y_2$.
We use the notation $\G = \G_X (Y)$, $\G_i= \G_X
(Y)_{<Y_i}$ for $i=1$ or $2$, and $\G_{12}$ for the intersection
$\G_1 \cap \G_2$.
Theorem \refT{Schlichting} can be applied to two inclusions of categories, $\G_{12} \to \G_1$ and $\G_{2} \to \G$.  The resulting homotopy fibrations are the rows in a commutative diagram of $K$-theory spectra
\[
\begin{CD}
K (\G_{12}) @>>> K (\G_1) @>>> K ({\G_1}/{\G_{12}}) \\
@VVV @VVV @VV{K(I)}V \\
K (\G_2) @>>> K (\G) @>>> K ({\G}/{\G_2})
\end{CD} \tag{$\natural$}
\]
The vertical maps are induced by exact inclusions, including $K(I)$ induced by the exact functor $I \colon {\G_1}/{\G_{12}} \to {\G}/{\G_2}$ which is itself induced from the exact inclusion $\G_1 \to \G$.

\begin{RemRef}{mbkdbdk}
	This is precisely the commutative diagram from Cardenas-Pedersen \cite[section 8]{mCeP:97} transported from bounded $K$-theory to fibred $G$-theory.  Cardenas and Pedersen use Karoubi quotients and the Karoubi fibrations in order to generate their diagram.  One of the crucial points in \cite{mCeP:97} is that the functor $I$ between the Karoubi quotients is an isomorphism of categories.
In fibred $G$-theory the situation is more complicated: $I$ is not necessarily full and, therefore, not an isomorphism of categories.  We will use the Approximation Theorem to prove that $K(I)$ is nevertheless an equivalence of spectra.
\end{RemRef}

\begin{PropRef}{THUIO}
$K(\wG') \simeq K({\G}/{\C})$.
\end{PropRef}

\begin{proof}
This follows from Lemma 2.3 in \cite{mS:03} as part of the proof of Theorem \refT{Schlichting} where $K(\wG')$ from Waldhausen's Fibration Theorem is identified with the Quillen $K$-theory spectrum $K ({\G}/{\C})$.
\end{proof}

\begin{LemRef}{CharWE}
If $f^{\subdot} \colon F^{\subdot} \to G^{\subdot}$ is a
degreewise admissible monomorphism with cokernel in $\C$  then
$f^{\subdot}$ is a weak equivalence in $\wG'$.
\end{LemRef}

\begin{proof}
The mapping cone $Cf^{\subdot}$ is quasi-isomorphic to the cokernel of
$f^{\subdot}$, by Lemma 11.6 of \cite{bK:96}, which is zero in
${\G}/{\C}$.
\end{proof}

The exact inclusion $I$ induces the exact functor $\wG'_1 \to \wG'$.

\begin{LemRef}{IdTarget}
The map $K  (\wG'_1) \to K (\wG')$ is a weak
equivalence.
\end{LemRef}

\begin{proof}
Applying the Approximation Theorem,
condition (1) is clear, so we need to check condition (2). Consider
\[
F^{\subdot} \colon \quad 0 \longrightarrow F^1 \xrightarrow{\
\phi_1\ } F^2 \xrightarrow{\ \phi_2\ }\ \dots\
\xrightarrow{\ \phi_{n-1}\ } F^n \longrightarrow 0
\]
in $\G_1$ and a chain map $g \colon F^{\subdot} \to G^{\subdot}$ for some complex
\[
G^{\subdot} \colon \quad 0 \longrightarrow G^1 \xrightarrow{\
\psi_1\ } G^2 \xrightarrow{\ \psi_2\ }\ \dots\
\xrightarrow{\ \psi_{n-1}\ } G^n \longrightarrow 0
\]
in $\G$.
Suppose all $F^i$ and $G^i$ are $(D,\Delta')$-lean/split and $(d,\delta)$-insular.
Also assume that there is a fixed number $r \ge 0$ and a monotone function $\rho \colon [0,+\infty) \to [0,+\infty)$ such that
$F^i (X,Y) \subset F^i ((X,C) [r,\rho_{x_0}])$
holds for all $0 \le i \le n$.
If the pair $(b,\theta)$ serves as bounded control data for all $\phi_i$, $\psi_i$, and $g_i$, we define the submodule
\[
F^{\prime{i}} = \mathcal{G}^{i} ((X,Y_1) [r+3ib,\rho_{x_0} + 3i\theta_{x_0}])
\]
and define $\xi_{i} \colon F^{\prime{i}} \to F^{\prime{i+1}}$ to be the restrictions of
$\psi_{i}$ to $F^{\prime{i}}$. This gives a chain subcomplex
$(F^{\prime{i}},\xi_{i})$ of $(G^{i},\psi_{i})$ in $\G$ with the inclusion $i \colon F^{\prime{i}} \to G^{i}$.  Notice that we have the induced chain map $\overline{g} \colon F^{\subdot} \to F^{\prime\subdot}$ in $\G_1$ so that $g = i I(\overline{g})$.

Once we establish that $C^{\subdot} = \coker (i)$ is in $\G_2$, $K(I)$ is a weak equivalence by Lemma \refL{CharWE}.

Since
\[
F^{\prime{i}} \subset {G}^{i} ((X,Y_1) [r+3ib +K,\rho_{x_0} + 3i\theta_{x_0} +k_{x_0}]),
\]
each $C^i$ is supported on
\begin{multline*}
(X,Y \setminus Y_1)[2D +2d - r - 3ib - K, 2\Delta'_{x_0} + 2\delta_{x_0} - \rho_{x_0} - 3i\theta_{x_0} - k_{x_0}] \\
\subset
(X, Y_2) [2D +2d, 2\Delta'_{x_0} + 2\delta_{x_0}],
\end{multline*}
cf.~Lemma \refL{SubSubLem1}.
So the complex $C^{\subdot}$ is indeed in $\G_2$.
\end{proof}

The excision theorems are best stated in terms of non-connective deloopings of the $K$-theory spectra.  Following Pedersen and Weibel we can use the same kind of diagram to first deloop $K$-theory and then reuse it to prove the excision theorem.

\medskip

Let $\mathbb{R}$, $\mathbb{R}^{\ge 0}$, and $\mathbb{R}^{\le 0}$
denote the metric spaces of the reals, the nonnegative reals, and
the nonpositive reals with the restriction of the usual metric on the
real line $\mathbb{R}$.  Then there is the following instance of
commutative diagram ($\natural$)

\[
\begin{CD}
{G}_X (Y) @>>> {G}_{X} (Y \times \mathbb{R}^{\ge 0}) @>>> K ({\G_1}/{\G_{12}}) \\
@VVV @VVV @VV{K(I)}V \\
{G}_X (Y \times \mathbb{R}^{\ge 0}) @>>> {G}_X (Y \times \mathbb{R}) @>>> K ({\G}/{\G_2})
\end{CD}
\]
\smallskip

We already know that $K(I)$ is an equivalence.

\begin{LemRef}{ESw}
The spectra ${G}_{X} (Y \times \mathbb{R}^{\ge 0})$ and ${G}_X (Y \times \mathbb{R}^{\le 0})$ are contractible.
\end{LemRef}

\begin{proof}
This follows from the fact that these controlled categories are
flasque, that is, the evident shift functor $T$ in the positive
(respectively negative)
direction along $\mathbb{R}^{\ge 0}$ (respectively $\mathbb{R}^{\le 0}$) interpreted in the obvious
way is an exact endofunctor, and there is a natural equivalence $1
\oplus \pm T \cong \pm T$. Contractibility follows from the
Additivity Theorem, cf.\ Pedersen--Weibel \cite{ePcW:85}.
\end{proof}

In view of Lemma \refL{IdTarget}, we
obtain a map ${G}_X (Y) \to \Omega {G}_{X} (Y \times \mathbb{R})$
which induces isomorphisms of $K$-groups
in positive dimensions.
Weak equivalences
\[
\Omega^k {G}_{X} (Y \times \mathbb{R}^k)
\longrightarrow \Omega ^{k+1}
{G}_{X} (Y \times \mathbb{R}^{k+1})
\]
are obtained by iterating this construction for $k \ge 2$.

\begin{DefRef}{RealExQ}
The \textit{nonconnective fibred bounded $G$-theory} over the pair $(X,Y)$ is the spectrum
\[
\GncX (Y) \overset{ \text{def} }{=} \hocolim{k>0}
\Omega^{k} {G}_{X} (Y \times \mathbb{R}^k).
\]
Since $\BLI (X,R)$ can be identified with $\B_X (\point)$, this definition also gives a nonconnective delooping of the
$G$-theory of $X$:
\[
\Gnc (X,R) = \hocolim{k > 0}
\Omega^{k} {G}_{X} (\mathbb{R}^k).
\]
\end{DefRef}

The subcategory $\G_{X} (Y \times \mathbb{R}^k)_{<C \times \mathbb{R}^k}$ is evidently a Serre subcategory of $\G_{X} (Y \times \mathbb{R}^k)$ for any choice of the subset $C \subset Y$.

\begin{DefRef}{RealExQr}
We define
\[
\GncX (Y)_{<C} \overset{ \text{def} }{=} \hocolim{k>0}
\Omega^{k} {G}_{X} (Y \times \mathbb{R}^k)_{<C \times \mathbb{R}^k}.
\]
We also define
\[
\GncX (Y)_{<C_1, C_2} \overset{ \text{def} }{=} \hocolim{k>0}
\Omega^{k} {G}_{X} (Y \times \mathbb{R}^k)_{<C_1 \times \mathbb{R}^k, \, C_2 \times \mathbb{R}^k}.
\]
\end{DefRef}

\begin{ThmRefName}{Exc2}{Fibrewise Bounded Excision, Version One}
Suppose $Y_1$ and $Y_2$ are subsets of a metric space $Y$, and $Y = Y_1 \cup Y_2$.
There is a homotopy pushout diagram of spectra
\[
\begin{CD}
\GncX (Y)_{<Y_1,Y_2} @>>> \GncX (Y)_{<Y_1} \\
@VVV @VVV \\
\GncX (Y)_{<Y_2} @>>> \GncX (Y)
\end{CD}
\]
where the maps of spectra are induced from the exact inclusions.
If $Y_1$ and $Y_2$ are mutually antithetic subsets of $Y$,
there is a homotopy pushout
\[
\begin{CD}
\GncX (Y_1 \cap Y_2) @>>> \GncX (Y_1) \\
@VVV @VVV \\
\GncX (Y_2) @>>> \GncX (Y)
\end{CD}
\]
\end{ThmRefName}

\begin{proof}
Let us write $S^k \G$ for $\G_{X} (Y \times \mathbb{R}^k)$
whenever $\G$ is the fibred bounded category for a pair $(X,Y)$.
If $C$ represents a family of coarsely equivalent subsets in a coarse covering $\mathcal{U}$ of $Y$, consider the
fibration
\[
G_X (C) \longrightarrow G_X (Y) \longrightarrow K ({\G}/{\C})
\]
from Theorem \refT{LocText}.  Notice that there is a map
$K ({\G}/{\C}) \to  \Omega K ({S\G}/{S\C})$
which is an equivalence in positive dimensions by the Five
Lemma.
Defining
\[
\GncX (Y,C) = \Knc ({\G}/{\C}) = \hocolim{k}
\Omega^{k} K ({S^k \G}/{S^k \C})
\]
gives an induced fibration
\[
\GncX (C) \longrightarrow \GncX (Y) \longrightarrow \GncX (Y,C).
\]
The theorem follows from the commutative diagram
\[
\begin{CD}
\GncX (Y)_{<Y_1,Y_2} @>>> \GncX (Y)_{<Y_1} @>>> \Knc ({\G_1}/{\G_{12}}) \\
@VVV @VVV @VVV \\
\GncX (Y)_{<Y_2} @>>> \GncX (Y) @>>> \Knc ({\G}/{\G_2})
\end{CD}
\]
and the fact that $\Knc  ({\G_1}/{\G_{12}})
\to \Knc ({\G}/{\G_2})$ is a weak equivalence.

From Lemma \refL{LemNeed1} one has a weak equivalence $\GncX (Y)_{<C} \simeq \GncX (C)$.  
If $C_1$ and $C_2$ are coarsely antithetic then the same construction shows that $\G_X (C_1 \cap C_2) \to \G_X (Y)_{<C_1, C_2}$ is onto the essential full image, and so there is a weak equivalence
$\GncX (Y)_{<C_1, C_2} \simeq \GncX (C_1 \cap C_2)$.
This allows to substitute the terms in the commutative diagram, giving the second statement.
\end{proof}

\SSecRef{Fibred coarse coverings}{FCCCCC}

To state the excision theorems properly in the coarse geometric setting, we develop the language of fibred coarse coverings.

Two subsets $A$, $B$ of $(X,Y)$ are called \textit{coarsely equivalent} if there is a set of enlargement data $(K,k)$ such that
$A \subset B[K,k_{x_0}]$ and $B \subset A[K,k_{x_0}]$. 
We will use the notation $A \, \| \, B$ for this equivalence relation.

A family of subsets $\mathcal{A}$ is called \textit{coarsely saturated} if it is maximal with respect to this equivalence relation.
Given a subset $A$, we denote by $\mathcal{S}(A)$ the smallest boundedly saturated family containing $A$.

A collection of subsets $\mathcal{U} = \{ U_i \}$ is a \textit{coarse covering} of $(X,Y)$ if $(X,Y) = \bigcup S_i$ for some $S_i \in \mathcal{S}(U_i)$.
Similarly, $\mathcal{U} = \{ \mathcal{A}_i \}$ is a \textit{coarse covering} by coarsely saturated families if for some (and therefore any) choice of subsets $A_i \in \mathcal{A}_i$, $\{ A_i \}$ is a coarse covering in the above sense.

We will say that a pair of subsets $A$, $B$ of $(X,Y)$ are \textit{coarsely antithetic} if
for any two sets of enlargement data $(D_1, d_1)$ and $(D_2, d_2)$ there exist enlargement data $(D, d)$ such that
\[
A[D_1, (d_1)_{x_0}] \cap B[D_2, (d_2)_{x_0}] \subset (A \cap B) [D, d_{x_0}].
\]
We will write $A \, \natural \, B$ to indicate that $A$ and $B$ are coarsely antithetic.

Given two subsets $A$ and $B$, we define
\[
\mathcal{S}(A,B) = \{ A' \cap B' \, \vert \, A' \in \mathcal{S}(A), B' \in \mathcal{S}(B), A' \, \natural \, B' \}.
\]
It is easy to see that $\mathcal{S}(A,B)$ is a coarsely saturated family.

\begin{PropRef}{UYNBVC}
$\mathcal{S}(A,B)$ is a coarsely saturated family.
\end{PropRef}

\begin{proof}
Suppose $A_1$, $A'_1$ and $A_2$, $A'_2$ are two coarsely antithetic pairs, and
$A_1 \subset A_2 [D_{12}, (d_{12})_{x_0}]$, $A'_1 \subset A'_2 [D'_{12}, (d'_{12})_{x_0}]$ for some $D_{12}$, $d_{12}$, $D'_{12}$, and $d'_{12}$.
Then
\[
A_1 \cap A'_1 \subset A_2 [D_{12}, (d_{12})_{x_0}] \cap A'_2 [D'_{12}, (d'_{12})_{x_0}] \subset (A_2 \cap A'_2) [D,d_{x_0}]
\]
for some $(D,d)$.
\end{proof}

There is the straightforward generalization to the case of a finite number of subsets of $(X,Y)$.
Similarly, we write $A_1 \, \natural \, \ldots \, \natural \, A_k$ if for arbitrary sets of data $(D_i,d_i)$ there is a set of enlargement data $(D,d)$ so that
\[
A_1 [D_1,(d_1)_{x_0}] \cap \ldots \cap A_k [D_k,(d_k)_{x_0}]
\subset (A_1 \cap \ldots \cap A_k) [D,d_{x_0}]
\]
and define
\[
\mathcal{S}(A_1, \ldots , A_k) = \{ A'_1 \cap \ldots \cap A'_k \, \vert \, A'_i \in \mathcal{S}(A_i), A_1 \, \natural \, \ldots \, \natural \, A_k \}.
\]
Identifying any coarsely saturated family $\mathcal{A}$ with $\mathcal{S}(A)$ for $A \in \mathcal{A}$, one has
the coarse saturated family
$\mathcal{S}(\mathcal{A}_1, \ldots , \mathcal{A}_k)$.
We will refer to $\mathcal{S}(\mathcal{A}_1, \ldots , \mathcal{A}_k)$ as the
\textit{coarse intersection} of $\mathcal{A}_1, \ldots , \mathcal{A}_k$.
A coarse covering $\mathcal{U}$ is \textit{closed under coarse intersections} if all coarse intersections $\mathcal{S}(\mathcal{A}_1, \ldots , \mathcal{A}_k)$ are nonempty and are contained in $\mathcal{U}$.
If $\mathcal{U}$ is a given coarse covering, the smallest coarse covering that is closed under coarse intersections and contains
$\mathcal{U}$ will be called the \textit{closure} of $\mathcal{U}$ under coarse intersections.

All of the terms introduced above have absolute analogues obtained by simply restricting to the case $Y = \point$.  So there are, in particular, finite coarse coverings of a single metric space.

\begin{PropRef}{UIVCXZSAW}
If $\mathcal{U}$ is a finite coarse antithetic covering of $Y$ then $(X,\mathcal{U})$ consisting of subsets $(X,U)$, $U \in \mathcal{U}$, is a coarse antithetic covering of $(X,Y)$.  If $\mathcal{U}$ is closed under coarse intersections, $(X,\mathcal{U})$ is closed under coarse intersections.
\end{PropRef}

\begin{proof}
Suppose $\mathcal{U} = \{ \mathcal{A}_i \}$ so that for $A_i \in \mathcal{A}_i$, $\{ A_i \}$ is a coarse covering of $Y$.  Then $\{ (X,A_i) \}$ is a covering of $(X,Y)$.  Suppose $\mathcal{U}$ is coarsely antithetic, so given numbers $d_1$, $d_2$ there is a number $d$ so that $A_i [d_1] \cap A_j [d_2] \subset (A_i \cap A_j) [d]$. If $d_1$, $d_2$ are non-decreasing functions, these values give a non-decreasing function $d$.  Now given enlargement data $(D_1, d_1)$ and $(D_2, d_2)$, we have 
\[
(X,A_i)[D_1, (d_1)_{x_0}] \cap (X,A_j)[D_2, (d_2)_{x_0}] \subset (X,A_i \cap A_j) [D, h],
\]
where $D$ can be any non-negative number, and $h$ is the function $h(x) = d (d_X (x_0, x) + D_1 +D_2)$. So $(X,\mathcal{U})$ is a coarsely antithetic covering.
A similar estimate gives the last statement.
\end{proof}

Suppose $\mathcal{U}$ is a finite coarse covering of $Y$ closed under coarse intersections.
We can define the homotopy pushout
\[
\mathcal{G}_X (Y; \mathcal{U}) = \hocolim{U \in \mathcal{U}} \GncX (Y)_{<U}.
\]

\begin{ThmRefName}{Exc2bb}{Fibrewise Bounded Excision, Version Two}
There is a weak equivalence
\[
\mathcal{G}_X (Y; \mathcal{U}) \simeq \GncX (Y).
\]
\end{ThmRefName}

\begin{proof}
Apply Theorem \refT{Exc2} inductively to the sets in $\mathcal{U}$.
\end{proof}

\SSecRef{Relative excision theorems}{REL}

Fibred $G$-theory has a useful relative version, and there are generalizations of the excision theorems to relative statements.  

\begin{DefRef}{GEPCWpr}
Let $Y' \in \mathcal{A}$ for a coarse covering $\mathcal{U}$ of $Y$.  Let $\G = \G_X (Y)_{< \mathcal{U}}$ and $\Y' = \G_X (Y)_{<Y'}$.  The category $\G_X (Y,Y')$ is the quotient category $\G/\Y'$.

It is now straightforward to define
\[
\GncX (Y,Y') {=} \hocolim{k>0}
\Omega^{k} {G}_{X} (Y \times \mathbb{R}^k, Y' \times \mathbb{R}^k),
\]
\[
\GncX (Y,Y')_{<C} {=} \hocolim{k>0}
\Omega^{k} {G}_{X} (Y \times \mathbb{R}^k, Y' \times \mathbb{R}^k)_{<C \times \mathbb{R}^k},
\]
and
\[
\GncX (Y,Y')_{<C_1, C_2} {=} \hocolim{k>0}
\Omega^{k} {G}_{X} (Y \times \mathbb{R}^k, Y' \times \mathbb{R}^k)_{<C_1 \times \mathbb{R}^k, \, C_2 \times \mathbb{R}^k}.
\]
\end{DefRef}

The theory developed in this section is spontaneously relativized to give the following excision theorem.

\begin{ThmRefName}{ExRel}{Relative Fibrewise Excision, Version One}
If $Y$ is the union of two subsets $U_1$ and $U_2$,
there is a homotopy pushout diagram of spectra
\[
\begin{CD}
\GncX (Y,Y')_{<U_1,U_2} @>>> \GncX (Y,Y')_{<U_1} \\
@VVV @VVV \\
\GncX (Y,Y')_{<U_2} @>>> \GncX (Y,Y')
\end{CD}
\]
where the maps of spectra are induced from the exact inclusions.
In fact, if $Y$ is the union of two mutually antithetic subsets $U_1$ and $U_2$, and $Y'$ is antithetic to both $U_1$ and $U_2$,
there is a homotopy pushout
\[
\begin{CD}
\GncX (U_1 \cap U_2,U_1 \cap U_2 \cap Y') @>>> \GncX (U_1,U_1 \cap Y') \\
@VVV @VVV \\
\GncX (U_2,U_2 \cap Y') @>>> \GncX (Y,Y')
\end{CD}
\]
\end{ThmRefName}

Finally, we want to state the relative excision theorem in the most familiar form.

\begin{PropRef}{ExcExcPrep}
Given a subset $U$ of $Y'$, there is a weak equivalence
\[
\GncX (Y,Y') \simeq \GncX (Y-U,Y'-U).
\]
\end{PropRef}

\begin{proof}
Consider the setup of Theorem \refT{Exc2} with $Y_1 = Y-U$ and $Y_2=Y'$, then Lemma \refL{IdTarget} shows that the map
\[
\frac{\G_X (Y)_{<(Y-U)}}{\G_X (Y)_{<(Y-U)} \cap \G_X (Y)_{<Y'}} \, \longrightarrow \, \frac{\G_X (Y)}{\G_X (Y)_{<Y'}}
\]
induces a weak equivalence on the level of $K$-theory.

Notice that, since $U$ is a subset of $Y'$, we have the interpretation
\[
{\G_X (Y)_{<(Y-U)} \cap \G_X (Y)_{<Y'}} = \G_X (Y)_{<(Y'-U)}.
\]
Now the maps of quotients
\[
\frac{\G_X (Y)}{\G_X (Y')} \, \longrightarrow  \, \frac{\G_X (Y)}{\G_X (Y)_{<Y'}}
\]
and
\[
\frac{\G_X (Y)_{<(Y-U)}}{\G_X (Y)_{<(Y'-U)}} \, \longleftarrow \, \frac{\G_X (Y-U)}{\G_X (Y'-U)}
\]
induced by fully faithful embeddings also induce weak equivalences.
Their composition gives the required equivalence.
\end{proof}

The relative theorem can be restated using coarse coverings in terms of the homotopy pushout
\[
\mathcal{G}_X (Y,Y'; \mathcal{U}) = \hocolim{U \in \mathcal{U}} \GncX (Y,Y')_{<U}.
\]

\begin{ThmRefName}{VBNCRT}{Relative Fibrewise Excision, Version Two}
	There is a weak equivalence
\[
\mathcal{G}_X (Y,Y'; \mathcal{U}) \simeq \GncX (Y,Y').
\]
\end{ThmRefName}

\begin{proof}
	Apply Theorem \refT{ExRel} inductively to the sets in $\mathcal{U}$.
\end{proof}

\SecRef{Conclusion}{JJJKK}

It is a familiar fact that $G$-theoretic approximations to the usually more meaningful $K$-theoretic invariants are easier to compute.  This paper confirms the pattern in the controlled algebra setting.  One approach to computing the $K$-theory leads one to consider $K$-theory with fibred control.  It is in this setting that the tools of bounded $K$-theory become insufficient for computation.  The paper \cite[section 5.2]{gCbG:18} contains an explicit example of failure of the standard localization tools in bounded $K$-theory based on Karoubi filtrations.  This paper, in contrast, uses a different technology of exact and Waldhausen categories.  The Fibrewise Excision Theorems from section 4 suffice to resolve the example from \cite{gCbG:18} and perform computations in more general geometric settings.  This material will appear in \cite{gCbG:19}, while the relationship between the $K$-theory of group rings for finitely generated groups and a $G$-theoretic analogue based on controlled $G$-theory is studied in \cite{gCbG:03,gCbG:15}. For the purpose of stating the results we restrict to regular coefficient rings $R$ of finite global dimension.  The conclusion is that the appropriate $G$-theory of the group ring is computable leveraging the results of this paper, while the Cartan comparison map from the $K$-theory is an equivalence for a remarkably large class of groups $\pi$ including all groups with finite $K(\pi,1)$ and finite decomposition complexity.

\end{document}